\documentclass[10pt]{amsart}
\usepackage{amssymb,amsbsy,amsmath,amsfonts,amssymb,amscd}
\usepackage{latexsym,euscript,exscale}

\title{Banach spaces without minimal subspaces}
\author {Valentin Ferenczi and Christian Rosendal}
\date {May 2007}
\newcommand {\ca} {{2^\N}}

\newcommand {\A}{\mathbb A}
\newcommand {\B}{\mathbb B}
\newcommand {\F}{\mathbb F}

\newcommand {\N}{\mathbb N}
\newcommand {\Q}{\mathbb Q}
\newcommand {\R}{\mathbb R}

\renewcommand{\leq}{\ensuremath{\leqslant}}
\renewcommand{\geq}{\ensuremath{\geqslant}}

\newcommand{\HH}{\mathbb H}
\newcommand{\PP}{\mathbb P}

\newcommand {\D}{\mathbb D}

\newcommand{\norm}[1]{\lVert#1\rVert}

\newcommand{\om}{\omega}
\newcommand{\eps}{\epsilon}

\newcommand{\iso}{\cong}

\newcommand{\tom} {\emptyset}

\newcommand{\saa}{\Rightarrow}
\newcommand{\equi}{\Leftrightarrow}

\newcommand{\til}{\rightarrow}

\newcommand{\Lim}[1]{\mathop{\longrightarrow}\limits_{#1}}

\newcommand {\del}{ \; \big| \;}

\newcommand {\ku} {\mathcal}

\newcommand{\inv}{^{-1}}

\newcommand {\e} {\exists}
\renewcommand {\a} {\forall}

\newcommand{\fed}{\boldsymbol}

\numberwithin{equation}{section}

\newtheorem{thm}{Theorem}[section]
\newtheorem{cor}[thm]{Corollary}
\newtheorem{lemme}[thm]{Lemma}
\newtheorem{prop} [thm] {Proposition}
\newtheorem{defi} [thm] {Definition}

\newtheorem{prob}[thm]{Problem}
\newtheorem{quest}[thm]{Question}
\newtheorem{ex}[thm]{Example}

\begin{document}
\thanks{The  second author was partially supported by NSF grant DMS 0556368 and by FAPESP}


\begin{abstract}
We prove three new dichotomies for Banach spaces \`a la W.T. Gowers'
dichotomies. The three dichotomies characterise respectively the spaces having
no minimal subspaces, having no subsequentially minimal basic sequences, and
having no subspaces crudely finitely representable in all of their subspaces.
We subsequently use these results to make progress on Gowers' program of
classifying Banach spaces by finding characteristic spaces present in every
space. Also, the results are used to embed any partial order of size $\aleph_1$
into the subspaces of any space without a minimal subspace ordered by
isomorphic embeddability.
\end{abstract}

 \maketitle

\tableofcontents

\section{Introduction}\label{intro}
In the paper \cite{g:dicho}, W.T. Gowers initiated a celebrated classification
theory for Banach spaces. Since the task of classifying all (even separable)
Banach spaces up to isomorphism is extremely complicated (just how complicated
is made precise in \cite{flr}), one may settle for a {\em loose classification
of Banach spaces up to subspaces}, that is, look for a list of classes of
Banach spaces such that:

(a) each class is {\em pure}, in the sense that if a space belongs to a class,
then every subspace belongs to the same class, or maybe, in the case when the
properties defining the class depend on a basis of the space, every block
subspace belongs to the same class,

(b) the classes are {\em inevitable}, i.e., every Banach space contains a subspace in one of the classes,

(c) any two classes in the list are disjoint,

(d) belonging to one class gives a lot of information about operators that may be defined on the space or on its subspaces.

\

We shall refer to this list as the {\em list of  inevitable classes of Gowers}.
Many classical problems are related to this classification program, as for
example the question whether every Banach space contains a copy of $c_0$ or
$\ell_p$, solved in the negative by B.S. Tsirelson in 1974 \cite{tsi}, or the
unconditional basic sequence problem, also solved negatively by Gowers and B.
Maurey in 1993 \cite{GM}.  Ultimately one would hope to establish such a list
so that any classical space appears in one of the classes, and so that
belonging to that class would yield most of the properties which are known for
that space. For example, any property, which is known for Tsirelson's space, is
also true for any of its block subspaces. So Tsirelson's space is a pure space,
and, as such, should appear in one of the classes with a reasonable amount of
its properties. Also, presumably the nicest among the classes would consist of
the spaces isomorphic to $c_0$ or $\ell_p$, $1 \leq p<\infty$.

\

After the discovery by Gowers and Maurey of the existence of {\em hereditarily
indecomposable} (or HI) spaces, i.e., spaces such that no subspace may be
written as the direct sum of infinite dimensional subspaces \cite{GM}, Gowers
proved that every Banach space contains either an HI subspace or a subspace
with an unconditional basis \cite{g:hi}. These were the first two examples of
inevitable classes. We shall call this dichotomy the {\em first dichotomy} of
Gowers. He then used his famous Ramsey or determinacy theorem \cite{g:dicho} to
refine the list by proving that any Banach space contains a subspace with a
basis such that either no two disjointly supported block subspaces are
isomorphic, or such that any two subspaces have further subspaces which are
isomorphic. He called the second property {\em quasi minimality}. This {\em
second dichotomy} divides the class of spaces with an unconditional basis into
two subclasses (up to passing to a subspace). Finally, recall that a space is
{\em minimal} if it embeds into any of its subspaces. A quasi minimal space
which does not contain a minimal subspace is called {\em strictly quasi
minimal}, so Gowers again divided the class of quasi minimal spaces into the
class of strictly quasi minimal spaces and the class of minimal spaces.

Obviously the division between minimal and strictly quasi-minimal spaces is not
a real dichotomy, since it does not provide any additional information. The
main result of this paper is to provide the missing dichotomy for minimality,
which we shall call the {\em third dichotomy}.

A first step in that direction was obtained by A.M. Pe\l czar, who showed that
any strictly quasi minimal space contains a further subspace with the
additional property of not containing any subsymmetric sequence \cite{anna}.
The first author proved that the same holds if one replaces subsymmetric
sequences by embedding homogeneous sequences (any subspace spanned by a
subsequence contains an isomorphic copy of the whole space) \cite{subsurfaces}.

A crucial step in the proofs of \cite{anna} and \cite{subsurfaces} is the
notion of {\em asymptoticity}. An asymptotic game of length $k$ in a space $E$
with a basis is a game where I plays integers $n_i$ and II plays block vectors
$x_i$ supported after $n_i$, and where the outcome is the length $k$ sequence
$(x_i)$. Asymptotic games have been  studied extensively and  the gap between
finite dimensional and infinite dimensional phenomena was usually bridged by
fixing a constant  and letting the length of the game tend to infinity. For
example, a basis is {\em asymptotic $\ell_p$}  if there exists $C$ such that
for any $k$, I has a winning strategy in the length $k$ asymptotic game so that
the outcome is $C$-equivalent to the unit vector basis of $\ell_p^k$.

In \cite{anna} it is necessary to consider asymptotic games of {\em infinite}
length, which are defined in an obvious manner. The outcome is then an infinite
block sequence. The proof of the theorem  in \cite{anna} is based on the
obvious fact that if a basic sequence $(e_i)$ is  subsymmetric and $(x_i)$ is a
block sequence of $(e_i)$, then II has a strategy in the infinite asymptotic
game in $E=[e_i]$ to ensure that the outcome is equivalent to $(x_i)$. In
\cite{subsurfaces} a similar fact for embedding homogeneous basic sequences is
obtained, but the proof is  more involved and a more general notion of
asymptoticity must be used. Namely, a {\em generalised asymptotic game} in a
space $E$ with a basis $(e_i)$ is a game where I plays integers $n_i$ and II
plays integers $m_i$ and vectors $x_i$ such that ${\rm supp}(x_i) \subseteq
[n_1,m_1]\cup\ldots\cup[n_i,m_i]$, and the outcome is the sequence $(x_i)$,
which may no longer be a block basis.

The second author analysed infinite asymptotic games in \cite{asymptotic} (a
previous study had also been undertaken by E. Odell and T. Schlumprecht in
\cite{os:trees}), showing that the most obvious necessary conditions are, in
fact, also sufficient for II to have a strategy to play inside a given set.
This was done through relating the existence of winning strategies to a
property of subspaces spanned by vectors of the basis with indices in some
intervals of integers. Now the methods of \cite{asymptotic} extend to the
setting of generalised asymptotic games and motivate the following definition.
A space $Y$ is {\em tight} in a basic sequence $(e_i)$ if there is a sequence
of successive intervals $I_0<I_1<I_2<\ldots$ of $\N$ such that for all infinite
subsets $A\subseteq \N$, we have
$$
Y\not\sqsubseteq [e_n\del n\notin \bigcup_{i\in A}I_i],
$$
where $Y\sqsubseteq X$ denotes that $Y$ embeds into $X$. In other words, any
embedding of $Y$ into $[e_i]$ has a ``large'' image with respect to
subsequences of the basis $(e_i)$ and cannot avoid an infinite number of the
subspaces  $[e_n]_{n\in I_i}$. We then define a {\em tight basis} as a basis
such that every subspace is tight in it and a {\em tight space} as a space with
a tight basis.

As we shall prove in Lemma \ref{strategic uniformisation}, using the techniques
of \cite{asymptotic}, essentially a block subspace $Y=[y_i]$ is not tight in
$(e_i)$, when II has a winning strategy in the generalised asymptotic game in
$[e_i]$ for producing a sequence equivalent to $(y_i)$. This relates the notion
of tight bases to the methods  of \cite{subsurfaces}, and by extending these
methods we prove the main result of this paper:

\begin{thm}[3rd dichotomy]\label{main}
Let $E$ be a Banach space without minimal subspaces. Then $E$ has a tight subspace.
\end{thm}

Theorem \ref{main} extends the theorems of \cite{anna,subsurfaces}, since it is
clear that a tight space cannot contain a subsymmetric or even embedding
homogeneous block-sequence. This dichotomy also provides an improvement to the
list of Gowers: a strictly quasi minimal space must contain a tight quasi
minimal subspace. Example \ref{ex} shows that this is a non-trivial refinement
of the unconditional and strictly quasi minimal class, and Corollary \ref{tsi}
states that Tsirelson's space is tight. Theorem \ref{main} also refines the
class of HI spaces in the list, i.e., every HI space contains a tight subspace,
although it is unknown whether the HI property for a space with a basis does
not already imply that the basis is tight.

\

Our actual examples of tight spaces turn out to satisfy  one of two stronger
forms of tightness. The first is called {\em tightness with constants}. A basis
$(e_n)$ is tight with constants when for every infinite dimensional space $Y$,
the sequence of successive intervals $I_0<I_1<\ldots$ of $\N$ witnessing the
tightness of $Y$ in $(e_n)$ may be chosen so that $Y \not\sqsubseteq_K [e_n
\del n \notin I_K]$ for each $K$. This is the case for Tsirelson's space.

The second kind of tightness is called {\em tightness by range}. Here the
range, ${\rm range} \ x$, of  a vector $x$ is the smallest interval of integers
containing its support, and the range  of a block subspace $[x_n]$ is
$\bigcup_n {\rm range} \ x_n$. A  basis $(e_n)$ is tight by range when for
every block subspace $Y=[y_n]$, the sequence of successive intervals
$I_0<I_1<\ldots$ of $\N$ witnessing the tightness of $Y$ in $(e_n)$ may be
defined by $I_k={\rm range}\ y_k$ for each $k$. This is equivalent to no two
block subspaces with disjoint ranges being comparable. In a companion paper
\cite{exemples}, we show that tightness by range is satisfied by an HI space
and also by a space with unconditional basis both constructed by Gowers.

It turns out that there  are natural dichotomies between each of these strong
forms of tightness and respective weak forms of minimality. For the first
notion, we define a space $X$ to be {\em locally minimal} if for some constant
$K$, $X$ is $K$-crudely finitely representable in any of its subspaces. Notice
that local minimality is easily incompatible with tightness with constants.
Using an equivalent form of Gowers' game, as defined by J. Bagaria and J.
L\'opez-Abad \cite{BL}, we prove:

\begin{thm}[5th dichotomy] \label{main3}
Any Banach space $E$ contains a subspace with a basis that is either tight
with constants or is locally minimal.
\end{thm}

The ideas involved in the notion of local minimality also make sense for block
representability, which allows us to connect these notions with asymptoticity
of basic sequences. Proving a simple dichotomy for when a space contains an
asymptotically $\ell_p$ subspace, we are led to the following dichotomy for
when a Banach space contains a copy of either $c_0$ or $\ell_p$.

\begin{thm}[The $c_0$ and $\ell_p$ dichotomy]\label{ellp}
Suppose $X$ is a Banach space not containing a copy of $c_0$ nor of $\ell_p$,
$1\leq p<\infty$. Then $X$ has a subspace $Y$ with a basis such that either
\begin{itemize}
  \item[(1)] $\a M \;\e n\;  \a U_1,\ldots, U_{2n}\subseteq Y\; \e u_i\in \ku
S_{U_i}$
$$
\Big(u_1<\ldots<u_{2n}\;\&\: (u_{2i-1})_{i=1}^n\not
\sim_M(u_{2i})_{i=1}^n\Big).$$
  \item[(2)] For all block bases $(z_n)$ of $Y=[y_n]$ there are
intervals $I_1 < I_2 < I_3 <\ldots$ such that $(z_n)_{n\in I_K}$ is not
$K$-equivalent to a block sequence of $(y_n)_{n\notin I_K}$.
\end{itemize}
Here, of course, the variables range over {\em infinite-dimensional} spaces.
\end{thm}
Property (1) indicates some lack of homogeneity and (2) some lack of
minimality. It is interesting to see which conditions the various examples of
Banach spaces not containing $c_0$ or $\ell_p$ satisfy; obviously, Tsirelson's
space and its dual satisfy (2) and indeed (2) is the only option for spaces
being asymptotic $\ell_p$. On the other hand, Schlumprecht's space $S$
satisfies (1).

\

There is also a dichotomy concerning tightness by range. This direction for
refining the list of inevitable classes of spaces was actually suggested by
Gowers in \cite{g:dicho}. P. Casazza proved that if a space $X$ has a shrinking
basis such that no block sequence is {\em even-odd} ( the odd subsequence is
equivalent to the even subsequence), then $X$ is not isomorphic to a proper
subspace, see \cite{g:hyperplanes}.  So any Banach space contains either a
subspace, which is not isomorphic to a proper subspace, or is saturated with
even-odd block sequences, and, in the second case, we may find a further
subspace in which Player II has a winning strategy to produce even-odd
sequences  in the game of Gowers associated to his Ramsey theorem. This fact
was observed by Gowers, but it was unclear to him what to deduce from the
property in the second case.

We answer this question by using Gowers' theorem to obtain a dichotomy which on
one side contains tightness by range, which is a slightly stronger property
than the Casazza property. On the other side, we define  a space $X$ with a
basis $(x_n)$ to be {\em subsequentially minimal} if every subspace of  $X$
contains an isomorphic copy of a subsequence of $(x_n)$. This last property is
satisfied by Tsirelson's space and will also be shown to be incompatible with
tightness by range.

\begin{thm}[4th dichotomy]\label{main2}
Any Banach space $E$ contains a
subspace with a basis that is either tight by range or is subsequentially minimal.
\end{thm}

It is easy to check that the second case in Theorem \ref{main2} may be improved
to the following hereditary property of a basis $(x_n)$, that we call {\em
sequential minimality}: every block sequence of $[x_n]$ has a further block
sequence $(y_n)$ such that every subspace of $[x_n]$ contains a copy of a
subsequence of $(y_n)$.

\

The five dichotomies and the interdependence of the properties involved can be
visualised in the following diagram.

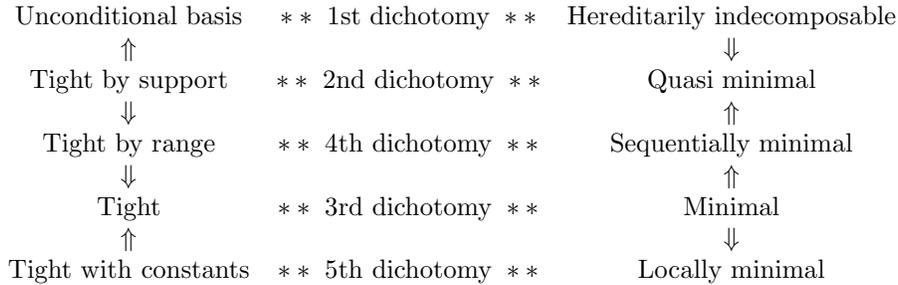
\[
\begin{tabular}{ccc}

Unconditional basis&$**\textrm{ 1st dichotomy }**$& Hereditarily indecomposable\\

$\Uparrow$&         &$\Downarrow$\\

Tight by support & $**\textrm{ 2nd dichotomy }**$ & Quasi minimal  \\

$\Downarrow$&&$\Uparrow$\\

Tight by range    &      $**\textrm{ 4th dichotomy }**$          & Sequentially minimal     \\

$\Downarrow$&&$\Uparrow$\\

Tight&                  $**\textrm{ 3rd dichotomy }**$              & Minimal\\

$\Uparrow$&         &$\Downarrow$\\

Tight with constants& $**\textrm{ 5th dichotomy }**$ & Locally minimal\\
\end{tabular}
\]

\

From a different point of view, coming from combinatorics and descriptive set
theory, Theorem \ref{main} also has important consequences for the isomorphic
classification of separable Banach spaces.  To explain this, suppose that $X$
is a Banach space and $SB_\infty(X)$ is the class of all infinite-dimensional
subspaces of $X$. Then the relation $\sqsubseteq$ of isomorphic embeddability
induces a partial order on the set of biembeddability classes of $SB_\infty(X)$
and we denote this partial order by $\PP(X)$. Many questions about the
isomorphic structure of $X$ translate directly into questions about the
structure of $\PP(X)$, e.g., $X$ has a minimal subspace if and only if $\PP(X)$
has a minimal element and $X$ is quasi minimal if and only if $\PP(X)$ is
downwards directed. In some sense, a space can be said to be pure in case the
complexity of $\PP(X)$ does not change by passing to subspaces and Gowers,
Problem 7.9 \cite{g:dicho}, motivated by this, asked for a classification of,
or at least strong structural information about, the partial orders $P$ for
which there is a Banach space $X$ saturated with subspaces $Y\subseteq X$ such
that $P\iso \PP(Y)$. A simple diagonalisation easily shows that such $P$ either
consist of a single point (corresponding to a minimal space) or are
uncountable, and, using methods of descriptive set theory and metamathematics,
this was successively improved in \cite{ergodic} and \cite{incomparable} to
either $|P|=1$ or $P$ having a continuum size antichain. Using a strengthening
of Theorem \ref{main}, we are now able to show that such $P$, for which
$|P|>1$, have an extremely complex structure by embedding any partial order of
size at most $\aleph_1$ into them.

For $A,B \subseteq \N$, we write $A \subseteq^* B$ to mean that $A \setminus B$
is finite.
\begin{thm}\label{posets}
Given a Banach space $X$, let $\PP(X)$ be the set of all
biembeddability classes of infinite-dimensional subspaces of $X$, partially ordered under isomorphic
embeddability. Let $P$ be a poset for which there exists a Banach space $X$ such that
$X$ is saturated with subspaces $Y$ such that $\PP(Y)\iso P$. Then either $|P|=1$, or
$\subseteq^*$ embeds into $P$. In the second case it follows that \begin{itemize}
  \item[(a)] any partial order of size at most $\aleph_1$ embeds into $P$, and
  \item[(b)] any closed partial order on a Polish space embeds into $P$.
\end{itemize}
\end{thm}
From the point of view of descriptive set theory, it is more natural to study
another problem, part of which was originally suggested to us by G. Godefroy
some time ago. Namely, the space $SB_\infty(X)$, for $X$ separable, can easily
be made into a standard Borel space using the Effros--Borel structure. In this
way, the relations of isomorphism, $\iso$, and isomorphic embeddability,
$\sqsubseteq$, become analytic relations on $SB_\infty(X)$ whose complexities
can be measured through the notion of Borel reducibility. We obtain Theorem
\ref{posets} as a consequence of some finer results formulated in this language
and that are of independent interest.

\

In Section \ref{gowers'dichotomies}, we put all the dichotomies together in
order to make progress on the loose classification mentioned above. In
connection with this, we shall also rely on work by A. Tcaciuc \cite{T}, who
proved a dichotomy for containing a strongly asymptotic $\ell_p$ basis, i.e., a
basis such that finite families of disjointly supported (but not necessarily
successive) normalised blocks supported ``far enough'' are uniformly equivalent
to the basis of $\ell_p^n$. Using just the first four dichotomies, in Theorem
\ref{gowersbis} we find 6 classes of inevitable spaces, 4 of which are known to
be non-empty, while if we use all 5 dichotomies plus Tcaciuc's, we find 19
classes. Out of these, 8 of them are known to be non-empty, though for 4 of the
examples, we will need the results of a companion paper \cite{exemples} where
these are constructed and investigated.

The resulting classification gives fairly detailed knowledge about the various
types of inevitable spaces, though much work remains to be done. In particular,
the new dichotomies explains some of the structural differences between the
wealth of new exotic spaces constructed in the wake of the seminal paper of
Gowers and Maurey \cite{GM}. It seems an interesting task to determine which of
the remaining 11 of the 19 cases are non-empty.

\section{Preliminaries}

\subsection{Notation, terminology, and conventions}
We shall in the following almost exclusively deal with infinite-dimensional
Banach spaces, so to avoid repeating this, we will always assume our spaces to
be infinite-dimensional. The spaces can also safely be assumed to be separable,
but this will play no role and is not assumed. Moreover, all spaces will be
assumed to be over the field of real numbers $\R$, though the results hold
without modification for complex spaces too.

Suppose $E$ is a Banach space with a normalised Schauder basis $(e_n)$. Then,
by a standard Skolem hull construction, there is a countable subfield ${\bf F}$
of $\R$ containing the rational numbers $\Q$ such that for any finite linear
combination
$$
\lambda_0e_0+\lambda_1e_1+\ldots+\lambda_ne_n
$$
with $\lambda_i\in {\bf F}$, we have
$\norm{\lambda_0e_0+\lambda_1e_1+\ldots+\lambda_ne_n}\in {\bf F}$. This means
that any ${\bf F}$-linear combination of $(e_n)$ can be normalised, while
remaining a ${\bf F}$-linear combination. Thus, as the set of $\Q$ and hence
also ${\bf F}$-linear combinations of $(e_n)$ are dense in $E$, also the set of
${\bf F}$-linear normalised combinations of $(e_n)$ are dense in the unit
sphere $\ku S_E$.

A {\em block vector} is a normalised finite linear combination
$x=\lambda_0e_0+\lambda_1e_1+\ldots+\lambda_ne_n$ where $\lambda_i\in {\bf F}$.
We insist on blocks being normalised and ${\bf F}$-linear and will be explicit
on the few occasions that we deal with non-normalised blocks. The restriction
to ${\bf F}$-linear combinations is no real loss of generality, but instead has
the effect that there are only countably many blocks. We denote by ${\bf D}$
the set of blocks. The {\em support}, ${\rm supp}\; x$, of a block
$x=\lambda_0e_0+\lambda_1e_1+\ldots+\lambda_ne_n$ is the set of $i\in \N$ such
that $\lambda_i\neq 0$ and the {\em range}, ${\rm range}\; x$, is the smallest
interval $I\subseteq \N$ containing ${\rm supp}\; x$.

A {\em block (sub)sequence}, {\em block basis}, or {\em blocking} of $(e_n)$ is
an infinite sequence $(x_n)$ of blocks such that ${\rm supp}\; x_n<{\rm supp}\;
x_{n+1}$ for all $n$ and a {\em block subspace} is the closed linear span of a
block sequence. Notice that if $X$ is a block subspace, then the associated
block sequence $(x_n)$ such that $X=[x_n]$ is uniquely defined up to the choice
of signs $\pm x_n$. So we shall sometimes confuse block sequences and block
subspaces. For two block subspaces $X=[x_n]$ and $Y=[y_n]$, write $Y\leq X$ if
$Y\subseteq X$, or, equivalently, $y_n\in {\rm span}(x_i)$ for all $n$. Also,
let $Y\leq ^* X$ if there is some $N$ such that $y_n\in {\rm span}(x_i)$ for
all $n\geq N$.

When we work with block subspaces of some basis $(e_n)$, we will assume that we
have chosen the same countable subfield ${\bf F}$ of $\R$ for all block
sequences $(x_n)$ of $(e_n)$, and hence a vector in $[x_n]$ is a block of
$(x_n)$ if and only if it is a block of $(e_n)$, so no ambiguity occurs. We
consider the set $bb(e_n)$ of block sequences of $(e_n)$ as a closed subset of
${\bf D}^\N$, where ${\bf D}$ is equipped with the discrete topology. In this
way, $bb(e_n)$ is a Polish, i.e., separable, completely metrisable space. If
$\Delta=(\delta_n)$ is a sequence of positive real numbers, which we denote by
$\Delta>0$, and $\A\subseteq bb(e_n)$, we designate by $\A_\Delta$ the set
$$
\A_\Delta=\{(y_n)\in bb(e_n)\del \e (x_n)\in bb(e_n)\; \a n\; \norm{x_n-y_n}<\delta_n\}.
$$

If $A$ is an infinite subset of $\N$, we denote by $[A]$ the space of infinite
subsets of $A$ with the topology inherited from $2^A$. Also, if $a\subseteq \N$
is finite,
$$
[a,A]=\{B\in [\N]\del a\subseteq B \subseteq a\cup ( A\cap  [\max a+1,\infty[) \}.
$$
We shall sometimes confuse infinite subsets of $\N$ with their increasing
enumeration. So if $A\subseteq \N$ is infinite, we denote by $A_n$ the $n+1$'st
element of $A$ in its increasing enumeration (we start counting at $0$).

A Banach space $X$ {\em embeds} into $Y$ if $X$ is isomorphic to a closed subspace of $Y$. Since we shall work with the embeddability relation as a mathematical object itself, we prefer to use the slightly non-standard notation $X\sqsubseteq Y$ to denote that $X$ embeds into $Y$.

Given two Banach spaces $X$ and $Y$, we say that $X$ is {\em crudely finitely
representable} in $Y$ if there is a constant $K$ such that for any
finite-dimensional subspace $F\subseteq X$ there is an embedding $T\colon F\til
Y$ with constant $K$, i.e., $\norm{T}\cdot\norm{T\inv}\leq K$.

Also, if $X=[x_n]$ and $Y=[y_n]$ are spaces with bases, we say that $X$ is {\em
crudely block finitely representable} in $Y$ if for some constant $K$ and all
$k$, there are (not necessarily normalised) blocks $z_0<\ldots<z_k$ of $(y_n)$
such that $(x_0,\ldots,x_k)\sim_K(z_0,\ldots,z_k)$.

Two Banach spaces are said to be {\em incomparable} if neither one
embeds into the other, and {\em totally incomparable} if no subspace of one is
isomorphic to a subspace of the other.

We shall at several occasions use non-trivial facts about the Tsirelson space
and its $p$-convexifications,
for which our reference is \cite{CS}, and also facts from descriptive set
theory that can all be found in \cite{kechris}. For classical facts in Banach
space theory we refer to \cite{LT}.

\subsection{Gowers' block sequence game}
A major ingredient in several of our proofs will be the following equivalent version of
Gowers' game due to J. Bagaria and J. L\'opez-Abad \cite{BL}.

Suppose $E=[e_n]$ is given. Player I and II alternate in choosing blocks $x_0<x_1<x_2<\ldots$ and $y_0<y_1<y_2<\ldots$ as follows:
Player I plays in the $k$'th round of the game a block $x_k$
such that $x_{k-1}<x_k$. In response to this, II either chooses to
pass, and thus play nothing in the $k$'th round, or plays a block
$y_i\in [x_{l+1},\ldots,x_k]$, where $l$ was the last round
in which II played a block.
$$
\begin{array}{cccccccccccccc}
{\bf I} & x_0 &\ldots  & x_{k_0} &  & x_{k_0+1} & \ldots&
x_{k_1}&\\
{\bf II}  & &  & &y_0\in [x_0,\ldots,x_{k_0}] &  &  &  & y_1\in[x_{k_0+1},\ldots,x_{k_1}]
\end{array}
$$
We thus see I as constructing a block sequence $(x_i)$, while II chooses a
block subsequence $(y_i)$. This block subsequence $(y_i)$ is then called the
{\em outcome} of the game. (Potentially the blocking could be finite, but the
winning condition can be made such that II loses unless it is infinite.) We now
have the following fundamental theorem of Gowers (though he only proves it for
real scalars, it is clear that his proof is valid for the field ${\bf F}$ too).

\begin{thm}[W.T. Gowers \cite{g:dicho}]
Suppose $(e_n)$ is a Schauder basis and $\A\subseteq bb(e_i)$ is an analytic
set such that any $(x_i)\in bb(e_i)$ has a block subsequence $(y_i)$ belonging
to $\A$, then for all $\Delta>0$, there is a block subsequence $(v_i)\in
bb(e_i)$ such that II has a strategy to play in $\A_\Delta$ if  I is restricted
to play blockings of $(v_i)$.
\end{thm}

\subsection{A trick and a lemma}\label{tricks}
We gather here a couple of facts that will be used repeatedly later on.

We shall at several occasions use coding with inevitable subsets of the unit
sphere of a Banach space, as was first done by L\'opez-Abad in \cite{lopez}. So
let us recall here the relevant facts and set up a framework for such codings.

Suppose $E$ is an infinite-dimensional  Banach space with a basis not
containing a copy of $c_0$. Then by the solution to the distortion problem by
Odell and Schlumprecht \cite{OS:distortion} there is a block subspace $[x_n]$
of $E$ and two closed subsets $F_0$ and $F_1$ of the unit sphere of $[x_n]$
such that ${\rm dist}(F_0,F_1)=\delta>0$ and such that for all block bases
$(y_n)$ of $(x_n)$ there are block vectors $v$ and $u$ of $(y_n)$ such that
$v\in F_0$ and $u\in F_1$. In this case we say that $F_0$ and $F_1$ are {\em
positively separated, inevitable, closed subsets of $\ku S_{[x_n]}$}.

We can now use the sets $F_0$ and $F_1$ to code infinite binary sequences,
i.e., elements of $2^\N$ in the following manner. If $(z_n)$ is a block
sequence of $(x_n)$ such that for all $n$, $z_n\in F_0\cup F_1$, we let
$\varphi((z_n))=\alpha\in \ca$ be defined by
\[
\alpha_n=\left\{
           \begin{array}{ll}
             0, & \hbox{if $z_n\in F_0$;} \\
             1, & \hbox{if $z_n\in F_1$.}
           \end{array}
         \right.
\]
Since the sets $F_0$ and $F_1$ are positively separated, this coding is fairly
rigid and can be extended to block sequences $(v_n)$ such that ${\rm
dist}(v_n,F_0\cup F_1)<\frac\delta2$ by letting $\varphi((v_n))=\beta\in\ca$ be
defined by
\[
\beta_n=\left\{
           \begin{array}{ll}
             0, & \hbox{if ${\rm dist}(v_n, F_0)<\frac\delta2$;} \\
             1, & \hbox{if ${\rm dist}(v_n, F_1)<\frac\delta2$.}
           \end{array}
         \right.
\]
In this way we have that if $(z_n)$ and $(v_n)$ are block sequences with
$z_n\in F_0\cup F_1$ and $\|v_n-z_n\|<\frac\delta2$ for all $n$, then
$\varphi((z_n))=\varphi((v_n))$.

One can now use elements of Cantor space $\ca$ to code other objects in various
ways. For example, let $\HH$ denote the set of finite non-empty sequences
$(q_0,q_1,\ldots,q_n)$ of rationals with $q_n\neq 0$. Then, as $\HH$ is
countable, we can enumerate it as $\vec h_0,\vec h_1,\ldots$. If now $(y_n)$
and $(v_n)$ are block sequences with
$\varphi((v_n))=0^{n_0}10^{n_1}10^{n_2}1\ldots$, then $(v_n)$ codes an infinite
sequence $\Psi((v_n),(y_n))=(u_n)$ of finite linear combinations of $(y_n)$ by
the following rule:
\[
u_k=q_0y_0+q_1y_1+\ldots+q_my_m,
\]
where $\vec h_{n_k}=(q_0,\ldots,q_m)$.

We should then notice three things about this type of coding:
\begin{itemize}
  \item [-] It is {\em inevitable}, i.e., for all block sequences $(y_n)$ of $(x_n)$
  and $\alpha\in \ca$, there is a block sequence $(v_n)$ of $(y_n)$ with
  $\varphi((v_n))=\alpha$.
  \item [-] It is {\em continuous}, i.e., to know an initial segment of
  $(u_n)=\Psi((v_n),(y_n))$, we only need to know initial segments of
  $(v_n)$ and of $(y_n)$.
  \item [-] It is {\em stable under small perturbations}. I.e., given $\eps>0$, we can
  find some $\Delta=(\delta_n)$ only depending on $\eps$ and the basis
  constant of $(x_n)$ with the following property. Assume that $(v_n)$ and
  $(y_n)$ are block bases of $(x_n)$ with $v_n\in F_0\cup F_1$ for all $n$
  and such that $\Psi((v_n),(y_n))=(u_n)$ is a block sequence of $(y_n)$
  with $\frac 12<\|u_n\|<2$. Then whenever $(v'_n)$ and $(y'_n)$ are other
  block sequences of $(x_n)$ with $\|v_n-v'_n\|<\frac\delta2$ and
  $\|y_n-y'_n\|<\delta_n$ for all $n$, the sequence
  $\Psi((v'_n),(y'_n))=(u'_n)$ will be a block sequence of $(y'_n)$ that is
  $1+\eps$-equivalent to $(u_n)$.
\end{itemize}

One can of course consider codings of other objects than sequences of vectors
and, depending on the coding, obtain similar continuity and stability
properties.

The inevitability of the coding is often best used in the following form.
\begin{itemize}
  \item [-] Suppose $\B$ is a set of pairs $((y_n),\alpha)$, where $(y_n)$ is a block
  sequence of $(x_n)$ and $\alpha\in \ca$, such that for all block
  sequences $(z_n)$ of $(x_n)$ there is a further block sequence $(y_n)$
  and an $\alpha$ such that $((y_n),\alpha)\in\B$. Then for all block
  sequences $(z_n)$ of $(x_n)$ there is a further block sequence $(y_n)$
  such that for all $n$, $y_{2n+1}\in F_0\cup F_1$ and
  $((y_{2n}),\varphi((y_{2n+1})))\in \B$.
\end{itemize}
To see this, let $(z_n)$ be given and notice that by the inevitability of the
coding there is a block sequence $(w_n)$ of $(z_n)$ such that $w_{3n+1}\in F_0$
and $w_{3n+2}\in F_1$. Pick now a block sequence $(v_n)$ of $(w_{3n})$ and an
$\alpha$ such that $((v_n),\alpha)\in \B$. Notice now that between $v_n$ and
$v_{n+1}$ there are block vectors $w_{3i_n+1}$ and $w_{3i_n+2}$ of $(z_n)$
belonging to $F_0$, respectively $F_1$. Thus, if we let $y_{2n}=v_n$ and set
$$
y_{2n+1}=\left\{
           \begin{array}{ll}
             w_{3i_n+1}, & \hbox{if $\alpha_n=0$;} \\
             w_{3i_n+2}, & \hbox{if $\alpha_n=1$.}
           \end{array}
         \right.
$$
then $((y_{2n}),\varphi((y_{2n+1})))\in \B$.

\begin{lemme}\label{block uniformisation}
Let $(x_n^0)\geq (x_n^1)\geq (x_n^2)\geq \ldots$ be a decreasing sequence of
block bases of a basic sequence $(x_n^0)$. Then there exists a block basis
$(y_n)$ of $(x_n^0)$ such that $(y_n)$ is $\sqrt K$-equivalent with a block
basis of $(x_n^K)$ for every $K\geq 1$.
\end{lemme}

\begin{proof}
Let $c(L)$ be a constant depending on the basis constant of $(x_n^0)$ such that
if two block bases differ in at most $L$ terms, then they are
$c(L)$-equivalent. Find now a sequence $L_1\leq L_2\leq \ldots$ of non-negative
integers tending to $+\infty$ such that $c(L_K)\leq \sqrt K$. We can now easily
construct an infinite block basis $(y_n)$ of $(x_n^0)$ such that for all $K\geq
1$ at most the first $L_K$ terms of $(y_n)$ are not blocks of
$(x_n^K)_{n=L_K+1}^\infty$. Then $(y_n)$ differs from a block basis of
$(x_n^K)$ in at most $L_K$ terms and hence is $\sqrt K$-equivalent with a block
basis of $(x_n^K)$.
\end{proof}

\section{Tightness}
\subsection{Tight bases}
The following definition is central to the rest of the paper.
\begin{defi}
Consider a Banach space $E$ with a basis $(e_n)$ and let $Y$ be an arbitrary
Banach space. We say that $Y$ is {\em tight in the basis} $(e_n)$ if there is a
sequence of successive non-empty intervals $I_0<I_1<I_2<\ldots$ of $\N$ such
that for all infinite subsets $A\subseteq \N$, we have
$$
Y\not\sqsubseteq [e_n\del n\notin \bigcup_{i\in A}I_i].
$$
In other words,  if $Y$ embeds into $[e_n]_{n\in B}$, then $B\subseteq \N$ intersects all but finitely many intervals $I_i$.

We say that $(e_n)$ is {\em tight} if every infinite-dimensional Banach space $Y$ is tight in $(e_n)$.

Finally, an infinite-dimensional Banach space $X$ is {\em tight} if it has a tight basis.
\end{defi}
Also, the following more analytical criterion will prove to be useful. For
simplicity, denote by $P_I$ the canonical projection onto $[e_n]_{n\in I}$.
\begin{lemme}\label{projection tightness}
Let $X$ be a Banach space, $(e_n)$ a basis for a space $E$, and $(I_n)$ finite
intervals such that $\min I_n\Lim{n\til \infty}\infty$ and for all infinite
$A\subseteq \N$,
$$
X\not\sqsubseteq[e_n]_{n\notin\bigcup_{k\in A}I_k}.
$$
Then whenever $T\colon X\til [e_n]$ is an embedding, we have $\liminf_k \norm{P_{I_k}T}>0$.
\end{lemme}

\begin{proof}
Suppose towards a contradiction that $T\colon X\til E$ is an embedding such
that for some infinite $A\subseteq \N$, $\lim_{\substack{k\til \infty\\k\in
A}}\norm{P_{I_k}T}=0$. Then, by passing to an infinite subset of $A$, we can
suppose that $\sum_{k\in A}\norm{P_{I_k}T}<\frac 12\norm{T\inv}\inv$ and that
the intervals $(I_n)_{n\in A}$ are disjoint. Thus, the sequence of operators
$(P_{I_k}T)_{k\in A}$ is absolutely summable and therefore the operator
$\sum_{k\in A}P_{I_k}T\colon X\til E$ exists and has norm
$<\frac12\norm{T\inv}\inv$.

But then for $x\in X$ we have
\begin{displaymath}
\norm{\sum_{k\in A}P_{I_k}Tx}\leq \norm{\sum_{k\in A}P_{I_k}T}\cdot\norm{x}\leq \frac
1{2\norm{T\inv}}\norm{x}\leq \frac 1{2\norm{T\inv}}\norm{T\inv}\cdot\norm{Tx}=\frac12\norm{Tx},
\end{displaymath}
and hence also
\begin{displaymath}
\norm{\big(T-\sum_{k\in A}P_{I_k}T\big)x}\geq \norm{Tx}-\norm{\sum_{k\in A}P_{I_k}Tx}\geq \norm{Tx}-\frac 12\norm{Tx}=\frac12\norm{Tx}.
\end{displaymath}
So $T-\sum_{k\in A}P_{I_k}T$ is still an embedding of $X$ into $E$. But this is
impossible as $T-\sum_{k\in A}P_{I_k}T$ takes values in
$[e_n]_{n\notin\bigcup_{k\in A}I_k}$. \end{proof}

\begin{prop}
A tight Banach space contains no minimal subspaces.
\end{prop}

\begin{proof}
Suppose $(e_n)$ is a tight basis for a space $E$ and let $Y$ be any subspace of
$E$. Pick a block subspace $X=[x_n]$ of $E$ that embeds into $Y$. Since $Y$ is
tight in $(e_n)$, we can find a sequence of intervals $(I_i)$ such that $Y$
does not embed into $[e_n]_{n\in B}$ whenever $B\subseteq \N$ is disjoint from
an infinite number of intervals $I_i$. By passing to a subsequence $(z_n)$ of
$(x_n)$, we obtain a space $Z=[z_n]$ that is a subspace of some $[e_n]_{n\in
B}$ where $B\subseteq \N$ is disjoint from an infinite number of intervals
$I_i$, and hence $Y$ does not embed into $Z$. Since $Z$ embeds into $Y$, this
shows that $Y$ is not minimal.
\end{proof}

The classical example of space without minimal subspaces is Tsirelson's space
$T$ and it is not too difficult to show that $T$ is tight. This will be proved
later on as a consequence of a more general result.

Any block sequence of a tight basis is easily seen to be tight. And also:
\begin{prop}
 If $E$ is a tight Banach space, then every shrinking basic sequence in $E$  is tight.
\end{prop}

\begin{proof}
Suppose $(e_n)$ is a tight basis for $E$ and $(f_n)$ is a shrinking basic
sequence in $E$. Let $Y$ be an arbitrary space and find intervals
$I_0<I_1<\ldots$ associated to $Y$ for $(e_n)$, i.e., for all infinite subsets
$A\subseteq \N$, we have $Y\not\sqsubseteq [e_n\del n\notin \bigcup_{i\in
A}I_i]$.

We notice that, since $(e_n)$ is a basis, we have for all $m$
$$
\norm{P_{I_k}|_{[f_i\del i\leq m]}}\Lim{k\til \infty}0, \eqno(1)
$$
and, since $(f_n)$ is shrinking and the $P_{I_k}$ have finite rank, we have for all $k$
$$
\norm{P_{I_k}|_{[f_i\del i> m]}}\Lim{m\til \infty}0. \eqno(2)
$$
Using alternately (1) and (2), we can construct integers $k_0<k_1<\ldots$ and intervals $J_0<J_1<\ldots$ such that
$$
\norm{P_{I_{k_n}}|_{[f_i\del i\notin J_n]}}<\frac 2{n+1}.
$$
To see this, suppose $k_{n-1}$ and $J_{n-1}$ have been defined and find some large $k_n>k_{n-1}$ such that
$$
\norm{P_{I_{k_n}}|_{[f_i\del i\leq \max J_{n-1}]}}\leq \frac 1{n+1}.
$$
Now, choose $m$ large enough that
$$
\norm{P_{I_{k_n}}|_{[f_i\del i> m]}}\leq\frac1{n+1},
$$
and set $J_n=[\max J_{n-1}+1,m]$. Then $\norm{P_{I_{k_n}}|_{[f_i\del i\notin
J_n]}}<\frac2{n+1}$. It follows that if $A\subseteq \N$ is infinite and
$T\colon Y\til [f_i]_{i\notin \bigcup_{n\in A}J_n}$ is an embedding, then
$\lim_{n\in A}\norm{P_{I_{k_n}}T}=0$, which contradicts Lemma \ref{projection
tightness}. So $(J_n)$ witnesses that $Y$ is tight in $(f_n)$.
\end{proof}

\begin{cor} If a tight Banach space $X$ is reflexive, then every basic sequence in $X$ is tight. \end{cor}

Notice that, since $c_0$ and $\ell_1$ are minimal, we have by the classical
theorem of James, that if $X$ is a tight Banach space with an unconditional
basis, then $X$ is reflexive and so every basic sequence in $X$ is tight.

\begin{ex}\label{ex}
The symmetrisation $S(T^{(p)})$  of the $p$-convexification $T^{(p)}$ of
Tsirel\-son's space, $1<p<+\infty$, does not contain a minimal subspace, yet it
is not tight.
\end{ex}
\begin{proof}
Since $S(T^{(p)})$ is saturated with isomorphic copies of subspaces of
$T^{(p)}$ and $T^{(p)}$ does not contain a minimal subspace, it follows that
$S(T^{(p)})$ does not have a minimal subspace. The canonical basis $(e_n)$ of
$S(T^{(p)})$ is symmetric, therefore $S(T^{(p)})$ is not tight in $(e_n)$ and
so $(e_n)$ is not tight. By reflexivity, no basis of $S(T^{(p)})$ is tight.
\end{proof}

\subsection{A generalised asymptotic game}
Suppose  $X=[x_n]$ and $Y=[y_n]$ are two Banach spaces with bases.
We define the game $H_{Y,X}$ with constant $C \geq 1$ between two
players I and II as follows: I will in each turn play a natural number $n_i$, while II will play
a not necessarily normalised block vector $u_i\in X$ and a natural number $m_i$ such that
$$
u_i\in X[n_0,m_0]+\ldots +X[n_i,m_i],
$$
where, for ease of notation,  we write $X[k,m]$ to denote $[x_n]_{k\leq n\leq m}$.
Diagramatically,
$$
\begin{array}{cccccccccccc}
{\bf I} & & & n_0 &  & n_1 &  & n_2 & & n_3& &\ldots \\
{\bf II} & & &  & u_0,m_0 &  & u_1, m_1 &  & u_2,m_2 & & u_3,m_3 &\ldots
\end{array}
$$
We say that the sequence $(u_i)_{i\in \N}$ is the {\em outcome} of the game and say that II wins the game if $(u_i)\sim_{C}(y_i)$.

\

For simplicity of notation, if $X=[x_n]$ is space with a basis, $Y$ a Banach
space,  $I_0<I_1<I_2<\ldots$  a sequence of non-empty intervals of $\N$ and $K$
is a constant, we write
$$
Y\sqsubseteq _K(X,I_i)
$$
if there is an infinite set $A\subseteq \N$ containing $0$ such that
$$
Y\sqsubseteq_K [x_n\del n\notin \bigcup_{i\in A}I_i],
$$
i.e., $Y$ embeds with constant $K$ into the subspace of $X$ spanned by $(x_n)_{n\notin \bigcup_{i\in A}I_i}$. Also, write
$$
Y\sqsubseteq (X,I_i)
$$
if there is an infinite set $A\subseteq \N$ such that $Y\sqsubseteq [x_n\del
n\notin \bigcup_{i\in A}I_i]$. Notice that in the latter case we can always
demand that $0\in A$ by perturbating the embedding with a finite rank operator.

It is clear that if $Y=[y_n]$ and II has a winning strategy in the game
$H_{Y,X}$ with constant $K$, then for any sequence of intervals $(I_i)$,
$Y \sqsubseteq_K(X,I_i)$.

Modulo the determinacy of open games, the next lemma shows that the converse holds up to a perturbation.
\begin{lemme}\label{strategic uniformisation}
Suppose $X=[x_n]$ is space with a basis and $K, \epsilon$  are positive
constants such that for all block bases $Y$ of $X$ there is a winning strategy
for I in the game $H_{Y,X}$ with constant $K+\epsilon$. Then there is a Borel
function $f\colon bb(X)\til [\N]$ such that for all $Y$ if
$I_j=[f(Y)_{2j},f(Y)_{2j+1}]$, then
$$
Y\not\sqsubseteq_K(X,I_j).
$$
\end{lemme}

\begin{proof}
Notice that the game $H_{Y,X}$ is open for player I and, in fact, if ${\bf
D}_{K+\epsilon}$ denotes the set of blocks $u$ with $\frac1{K+\epsilon}\leq
\norm{u}\leq K+\epsilon$, then the set
\begin{displaymath}\begin{split}
\A=\{&(Y,\vec p)\in bb(X)\times (\N\times {\bf D}_{K+\epsilon}\times \N)^\N\del
\textrm{ either $\vec p$ is a legal run of the game } H_{Y,X}\\
&\textrm{with constant $K+\epsilon$ in which I wins or $\vec p$ is not a legal run of the game } H_{Y,X}\}
\end{split}\end{displaymath}
is Borel and has open sections $\A_Y=\{\vec p\in (\N\times {\bf
D}_{K+\epsilon}\times \N)^\N\del (Y,\vec p)\in \A\}$. Also, since there are no
rules for the play of I in $H_{Y,X}$, $\A_Y$ really corresponds to the winning
plays for I in $H_{Y,X}$ with constant $K+\epsilon$. By assumption, I has a
winning strategy to play in $\A_Y$ for all $Y$, and so by the theorem on
strategic uniformisation (see (35.32) in \cite{kechris}), there is a Borel
function $\sigma\colon Y\mapsto \sigma_Y$ that to each $Y$ associates a winning
strategy for I in the game $H_{Y, X}$ with constant $K+\epsilon$.

Now let $\Delta=(\delta_n)$ be a sequence of positive reals such that for all
$2KC$-basic sequences of  blocks $(w_n)$ of $X$ with $\frac1{K}\leq
\norm{w_n}\leq K$ (where $C$ is the basis constant of $X$) and sequences of
vectors $(u_n)$, if for all $n$, $\norm{w_n-u_n}<\delta_n$, then
$(w_n)\sim_{\sqrt{1+\epsilon/K}}(u_n)$. We also choose sets $\D_n$ of finite
(not necessarily normalised) blocks with the following properties:
\begin{itemize}
  \item [-] for each finite $d\subseteq \N$, the number of vectors $u\in \D_n$ such that ${\rm supp}\; u=d$ is finite,
  \item [-] for all blocks vectors $w$ with $\frac1{K}\leq \norm{w}\leq K$, there is
  some $u\in \D_n$ with ${\rm supp} \;w={\rm supp}\; u$ such that
  $\norm{w-u}<\delta_n$.
\end{itemize}
This is possible since the $K$-ball in $[x_i]_{i\in d}$ is totally bounded for
all finite $d\subseteq\N$. So for all $2KC$-basic sequences $(w_n)$ of blocks
with $\frac1K\leq \norm{w_n}\leq K$, there is some $(u_n)\in \prod_n\D_n$ such
that ${\rm supp}\; w_n={\rm supp}\; u_n$ and $\norm{w_n-u_n}<\delta_n$ for all
$n$, whence $(w_n)\sim_{\sqrt{1+\epsilon/K}}(u_n)$.

Suppose now that $Y=[y_n]$ is given.
For each  $p=(n_0,u_0,m_0,\ldots,n_i,u_i,m_i)$, where $u_j\in \D_j$ for all $j$ and
$$
\begin{array}{cccccccccccc}
{\bf I} & & & n_0 &                 & n_1 &                 & \ldots &  n_i& & \\
{\bf II} & & &       & u_0,m_0 &        & u_1, m_1 & \ldots & & & u_i,m_i&
\end{array}
$$
is a legal position in the game $H_{Y,X}$ in which I has played according to
$\sigma_Y$, we write $p<k$ if $n_j,u_j,m_j<k$ for all $j\leq i$. Notice that
for all $k$ there are only finitely many such $p$ with $p<k$, so we can define
$$
\alpha(k)=\max(k,\max\{\sigma_Y(p)\del p<k\})
$$
and set $I_k=[k,\alpha(k)]$.
Clearly, the sequence $(I_k)$ can be computed in a
Borel fashion from $Y$. The $I_k$ are not necessarily successive, but their
minimal elements tend to $\infty$, so to prove the lemma it is enough to show
that $Y$ does not $K$-embed into $[x_n]$ avoiding an infinite number of $I_k$
including $I_0$.

Suppose now for a contradiction that $A\subseteq \N$ is infinite, $0\in A$ and
$y_i\mapsto w_i$ is a $K$-embedding into $[x_n\del n\notin \bigcup_{k\notin
  A}I_k]$. By perturbing the embedding slightly, we can suppose that the $w_i$
are blocks such that $\frac 1K\leq \norm{w_i}\leq K$ and we still have a
$K\sqrt{1+\epsilon/K}$-embedding. Using the defining properties of $\D_i$, we
find $u_i\in \D_i$ such that $\norm{w_i-u_i}<\delta_i$ and ${\rm
  supp}\;w_i={\rm supp}\; u_i$ for all $i$, whereby
$(u_i)\sim_{\sqrt{1+\epsilon/K}}(w_i)\sim_{K\sqrt{1+\epsilon/K}}(y_i)$, and
therefore $(u_i) \sim_{K+\epsilon} (y_i)$.

We now proceed to define natural numbers $n_i$, $m_i$, and $a_i\in A$ such that
for $p_i=(n_0,u_0,m_0,\ldots, n_i,u_i,m_i)$, we have
\begin{itemize}
  \item [(i)] $a_0=0$ and $[0,n_0[\subseteq I_{a_0}$,
  \item [(ii)] $m_i=a_{i+1}-1$,
  \item [(iii)] $p_i$ is a legal position in $H_{Y,X}$ in which I has played according to $\sigma_Y$,
  \item [(iv)] $]m_i,n_{i+1}[\subseteq I_{a_{i+1}}$.
\end{itemize}
Let $a_0=0$ and $n_0=\sigma_Y(\emptyset)=\alpha(0)$, whence
$I_{a_0}=[0,\alpha(0)]=[0,n_0]$. Find $a_1$ such that $n_0,u_0,a_0<a_1$ and set
$m_0=a_1-1$. Then $p_0=(n_0,u_0,m_0)$ is a legal position in $H_{Y,X}$ in which
I has played according to $\sigma_Y$, $p_0<a_1$, so
$n_1=\sigma_Y(n_0,u_0,m_0)\leq \alpha(a_1)$, and therefore $]m_0,n_1[\subseteq
I_{a_1}=[a_1,\alpha(a_1)]$.

Now suppose by induction that $n_0,\ldots, n_i$ and $a_0,\ldots, a_i$ have been
defined. Since $[0,n_0[\subseteq I_{a_0}$ and $]m_j,n_{j+1}[\subseteq
I_{a_{j+1}}$ for all $j<i$, we have
$$
u_{i}\in X[n_0,m_0]+\ldots + X[n_{i-1},m_{i-1}]+ X[n_i,\infty[.
$$
Find some $a_{i+1}$  greater than all of $n_0,\ldots, n_i$, $u_0,\ldots, u_i$,
$a_0,\ldots, a_i$ and let $m_i=a_{i+1}-1$. Then
$$
u_{i}\in X[n_0,m_0]+\ldots + X[n_{i-1},m_{i-1}]+ X[n_i,m_i]
$$
and $p_i=(n_0,u_0,m_0,\ldots, n_i,u_i,m_i)$ is a legal position played
according to $\sigma_Y$. Since $p_i<a_{i+1}$ also
$$
n_{i+1}=\sigma_Y(n_0,u_0,m_0,\ldots, n_i,u_i,m_i)\leq \alpha(a_{i+1}).
$$
Thus $]m_i,n_{i+1}[\subseteq I_{a_{i+1}}=[a_{i+1},\alpha(a_{i+1})]$.

Now since $p_0\subseteq p_1\subseteq p_2\subseteq\ldots$, we can let $\vec
p=\bigcup_ip_i$ and see that $\vec p$ is a run of the game in which I followed
the strategy $\sigma_Y$ and II has played $(u_i)$. Since $\sigma_Y$ is winning
for I, this implies that $(u_i)\not\sim_{K+\epsilon}(y_i)$ contradicting our
assumption.
\end{proof}

\begin{lemme}\label{interval diagonalisation}
Suppose $X=[x_n]$ is a space with a basis and $Y$ is a space such that for all
constants $K$ there are intervals $I_0^{(K)}<I_1^{(K)}<I^{(K)}_2<\ldots$ such
that $Y\not\sqsubseteq_K(X,I^{(K)}_j)$. Then there are intervals
$J_0<J_1<J_2<\ldots$ such that $Y\not\sqsubseteq (X,J_j)$. Moreover, the
intervals $(J_j)$ can be computed in a  Borel manner from $(I_i^{(K)})_{i,K}$.
\end{lemme}

\begin{proof}
By induction we can  construct intervals $J_0<J_1<J_2<\ldots$ such that $J_n$
contains one interval from each of $(I^{(1)}_i), \ldots , (I^{(n)}_i)$ and if
$M=\min J_n-1$ and $K=\lceil n\cdot c(M)\rceil$, then $\max J_n>\max
I^{(K)}_0+M$, where $c(M)$ is a constant such that if two subsequences of
$(x_n)$ differ in at most $M$ terms then they are $c(M)$ equivalent. It then
follows that if $A\subseteq \N$ is infinite, then
$$
Y\not\sqsubseteq [x_n]_{n\notin \bigcup_{i\in A}J_i}.
$$
To see this, suppose towards a contradiction that $A\subseteq \N$ is infinite and that for some integer $N$,
$$
Y\sqsubseteq_N  [x_n]_{n\notin \bigcup_{i\in A}J_i}.
$$
Pick then $a\in A$ such that $a\geq N$ and set $M=\min J_a-1$ and $K=\lceil
a\cdot c(M)\rceil$. Define an isomorphic embedding $T$ from
$$
[x_n\del n\notin \bigcup_{i\in A}J_i]
$$
into
$$
[x_n\del n\notin \bigcup_{i\in A}J_i\;\&\; n>\max J_a]+[x_n\del \max I_0^{(K)}<n\leq \max J_a ]
$$
by setting
$$
T(x_n)=\left\{
                                \begin{array}{ll}
                                  x_n, & \hbox{if $n>\max J_a$;} \\
                                  x_{\max I_0^{(K)}+n+1}, & \hbox{if $n\leq M$.}
                                \end{array}
                              \right.
$$
This is possible since $\max J_a>\max I^{(K)}_0+M$. Also, since $T$ only
changes at most $M$ vectors from $(x_n)$, it is a $c(M)$ embedding. Therefore,
by composing with $T$ and using that $N\cdot c(M)\leq a\cdot c(M)\leq K$, we
see that
$$
Y\sqsubseteq_K [x_n\del n\notin \bigcup_{i\in A}J_i\;\&\; n>\max J_a]+[x_n\del \max I_0^{(K)}<n\leq \max J_a].
$$
In particular, as almost all $J_i$ contain an interval $I_l^{(K)}$, we can find
and infinite set $B\subseteq \N$ containing $0$ such that
$$
Y\sqsubseteq_K [x_n\del n\notin \bigcup_{i\in B}I^{(K)}_i],
$$
which is a contradiction.
\end{proof}

\begin{lemme}\label{borel tight}
Let $E=[e_n]$ be given and suppose that for all block subspaces $Z\leq E$ and
constants $C$ there is a block subspace $X\leq Z$ such that for all block
subspaces $Y\leq X$, I has a winning strategy in the game $H_{Y,X}$ with
constant $C$. Then there is a block subspace $X\leq E$ and a Borel function
$f\colon bb(X)\til [\N]$ such that for all normalised block bases $Y\leq X$, if
we set $I_j=[f(Y)_{2j},f(Y)_{2j+1}]$, then
$$
Y\not\sqsubseteq (X,I_j).
$$
\end{lemme}

\begin{proof}
Using the hypothesis inductively together with Lemma \ref{strategic
uniformisation}, we can construct a sequence $X_0\geq X_1\geq X_2\geq \ldots$
of block subspaces $X_K$ and corresponding Borel functions $f_K\colon
bb(X_K)\til [\N]$ such that for all $V\leq X_K$ if
$I_j=[f_K(V)_{2j},f_K(V)_{2j+1}]$, then $V\not\sqsubseteq_{K^2} (X_K,I_j)$.

Pick by Lemma \ref{block uniformisation} some block $X_\infty$ of $X_0$ that is
$\sqrt K$-equivalent with a block sequence $Z_K$ of $X_K$ for every $K\geq 1$.
Then for any block sequence $Y$ of $X_\infty$ and any $K\geq 1$ there is some
block sequence $V\leq Z_K\leq X_K$ such that $Y$ is $\sqrt K$-equivalent with
$V$. Let $(I_j)$ be the intervals given by $f_K(V)$ so that
$V\not\sqsubseteq_{K^2}(X_K,I_j)$. We can then in a Borel way from $(I_j)$
construct intervals $(J_j)$ such that $V\not\sqsubseteq_{K^2}(Z_K,J_j)$ and
therefore also $Y\not\sqsubseteq_{K}(X_\infty,J_j)$.

This means that there are Borel functions $g_K\colon bb(X_\infty)\til [\N]$
such that for all $Y\leq X_\infty$ if $J^K_j(Y)=[g_K(Y)_{2j},g_K(Y)_{2j+1}]$,
then $Y\not\sqsubseteq_K (X_\infty, J^K_j(Y))$. Using Lemma \ref{interval
diagonalisation} we can now in a Borel manner in $Y$ define intervals
$L^Y_0<L^Y_1<\ldots$ such that
$$
Y\not\sqsubseteq (X_\infty,L^Y_j).
$$
Letting $f\colon bb(X_\infty)\til [\N]$ be the Borel function corresponding to $Y\mapsto (L^Y_j)$, we have our result.
\end{proof}

As will be clear in Section \ref{chains and posets} it can be useful to have a
version of tightness that not only assures us that certain intervals exist, but
also tells us how to obtain these. Thus, we  call a basis $(e_n)$ {\em
continuously tight} if there is a continuous function $f\colon bb(e_n)\til
[\N]$ such that for all normalised block bases $X$, if we set
$I_j=[f(X)_{2j},f(X)_{2j+1}]$, then
$$
X\not\sqsubseteq ([e_n],I_j),
$$
i.e., $X$ does not embed into $[e_n]$ avoiding an infinite number of the intervals $I_j$.

We shall now improve Lemma \ref{borel tight} to conclude continuous tightness from its hypothesis.
\begin{lemme}\label{cont tight}
Let $E=[e_n]$ be  given and suppose that for all block subspaces $Z\leq E$ and
constants $C$ there is a block subspace $X\leq Z$ such that for all block
subspaces $Y\leq X$, I has a winning strategy in the game $H_{Y,X}$ with
constant $C$. Then there is a continuously tight block subspace $X\leq E$.
\end{lemme}

\begin{proof}
We observe that $E$ does not contain a copy of $c_0$. Indeed if $Z$ is a block
subspace of $E$ spanned by a block sequence which is $C$-equivalent to the unit
vector basis of $c_0$, then for any $Y \leq X \leq Z$, II has a winning
strategy in the game $H_{Y,X}$ with constant $C^2$. We shall then use codings
with inevitable subsets. So find first a block subspace $Z$ of $E$ such that
there are inevitable, positively separated, closed subsets $F_0$ and $F_1$ of
$\ku S_Z$. By Lemma \ref{borel tight}, we can find a further block subspace $V$
of $Z$ and and a Borel function $g\colon bb(V)\til [\N]$ such that for all
$Y\leq V$, if  $I_j=[g(Y)_{2j},g(Y)_{2j+1}]$, then $Y\not\sqsubseteq (V,I_j)$.
Define the set
\[
\begin{split}
\A=\big\{&(y_n)\in bb(V)\del y_{2n}\in F_0\equi n\notin g((y_{2n+1})) \textrm{ and } y_{2n}\in F_1\equi n\in g((y_{2n+1})\big\}.
\end{split}
\]
Obviously, $\A$ is Borel and, using inevitability, one can check that any block
basis of $V$ contains a further block basis in $\A$. Thus, by  Gowers'
Determinacy Theorem, we have that for all $\Delta>0$ there is a block sequence
$X$ of $V$ such that II has a strategy to play into $\A_\Delta$ when I plays
block subspaces of $X$. Choosing $\Delta>0$ sufficiently small, this easily
implies that for some block basis $X$ of $E$, there is a continuous function
$h\colon bb(X)\til bb(X)\times [\N]$ that to each $W\leq X$ associates a pair
$\big(Y,(I_n)\big)$ consisting of a block sequence $Y$ of $W$ and a sequence of
intervals $(I_n)$ such that $Y\not\sqsubseteq (V,I_j)$. Notice now that
continuously in the sequence $(I_j)$, we can construct intervals $(J_j)$ such
that $Y\not\sqsubseteq (X,J_j)$ and hence also $W\not\sqsubseteq (X,J_j)$. So
the continuous function $f\colon bb(X)\til [\N]$ corresponding to $W\mapsto
(J_j)$ witnesses the continuous tightness of $X$.
\end{proof}
We shall need the following consequence of continuous tightness in Section \ref{chains and posets}.
\begin{lemme}\label{ramsey tight}
Suppose $(e_n)$ is continuously tight. Then there is a continuous function
$f\colon [\N]\til [\N]$ such that for all $A, B\in [\N]$,  if $B$ is disjoint
from an infinite number of intervals $[f(A)_{2i},f(A)_{2i+1}]$, then
$[e_n]_{n\in A}$ does not embed into $[e_n]_{n\in B}$.
\end{lemme}

\begin{proof}
It is enough to notice that the function $h\colon[\N]\til bb(e_n)$ given by
$h(A)=(e_n)_{n\in A}$ is continuous. So when composed with the function
witnessing continuous tightness we have the result.
\end{proof}

\subsection{A game for minimality}
For $L$ and $M$  two block subspaces of $E$, define the infinite game $G_{L,M}$
with constant $C\geq 1$ between two players as follows. In each round I chooses
a subspace $E_i\subseteq L$ spanned by a finite block sequence of $L$, a
normalised block vector $u_i\in E_0+\ldots+E_i$, and an integer $m_i$. In the
first round II plays an integer $n_0$, and in all subsequent rounds II plays a
subspace $F_i$ spanned by a finite block sequence of $M$, a (not necessarily
normalised) block vector $v_i\in F_0+\ldots+F_i$ and an integer $n_{i+1}$.
Moreover, we demand that $n_i\leq E_i$ and $m_i\leq F_i$.

Diagramatically,
\[\begin{array}{ccccccc}
{\bf I}   &       &n_0\leq E_0 \subseteq  L &      &  n_1\leq E_1\subseteq L&     & \ldots     \\
              &       &u_0\in E_0,        m_0         &       & u_1\in E_0+E_1,  m_1          &       &      \\
              &      &               &      &            &     &            \\
{\bf II}   &n_0      &        & m_0\leq F_0 \subseteq M&       &m_1\leq F_1 \subseteq M & \ldots   \\
                       &      &           &v_0 \in F_0, n_1 &       &v_1 \in F_0+F_1, n_2&     \\
 \end{array}\]
The {\em outcome} of the game is the pair of infinite sequences $(u_i)$ and
$(v_i)$ and we say that II wins the game if $(u_i)\sim_C(v_i)$.

\begin{lemme}\label{games}
Suppose that $X$ and $Y$ are block subspaces of $E$ and that player II has a
winning strategy in the game $H_{Y,X}$ with constant $C$. Then II has a winning
strategy in the game $G_{Y,X}$ with constant $C$.
\end{lemme}

\begin{proof}
We shall in fact prove that II has a winning strategy in a game that is
obviously harder for her to win. Namely, we shall suppose that II always plays
$n_i=0$, which obviously puts less restrictions on the play of I. Moreover, we
do not require I to play the finite-dimensional spaces $E_i$, which therefore
also puts fewer restrictions on I in subsequent rounds. Therefore, we shall
suppress all mention of $E_i$ and $n_i$ and only require that the $u_i$ are
block vectors in $Y$.

While playing the game $G_{Y,X}$, II will keep track of an auxiliary play of
the game $H_{Y,X}$ in the following way. In the game $G_{Y,X}$ we have the
following play
\[\begin{array}{ccccccc}
{\bf I}        &u_0\in  Y, m_0&      &  u_1\in Y, m_1&     & \ldots     \\
                   &               &      &            &     &            \\
{\bf II}       &        & m_0\leq F_0 \subseteq X&       &m_1\leq F_1 \subseteq X & \ldots   \\
                   &           &v_0 \in F_0 &       &v_1 \in F_0+F_1&     \\
 \end{array}\]
We write each vector $u_i=\sum_{j=0}^{k_i}\lambda_j^iy_j$ and may for
simplicity of notation assume that $k_i<k_{i+1}$. The auxiliary run of
$H_{Y,X}$ that II will keep track of is as follows, where II plays according to
her winning strategy for $H_{Y,X}$.
$$
\begin{array}{cccccccccccc}
{\bf I} & m_0 && \ldots & m_0 &  &m_1 & &\ldots& m_1& &\ldots \\
{\bf II}  &  & w_0,p_0 & \ldots& & w_{k_0}, p_{k_0} &  & w_{k_0+1},p_{k_0+1} &\ldots && w_{k_1},p_{k_1} &\ldots
\end{array}
$$
To compute the $v_i$ and $F_i$ in the game $G_{Y,X}$, II will refer to the play of $H_{Y,X}$ and set
$$
v_i=\sum_{j=0}^{k_i}\lambda_j^iw_j,
$$
and let
$$
F_i=X[m_i,\max\{p_{k_{i-1}+1},\ldots,p_{k_i}\}].
$$
It is not difficult to see that $m_i\leq F_i\subseteq X$, $v_i\in F_0+\ldots
+F_i$, and that the $F_i$ and $v_i$ only depends on $u_0,\ldots, u_i$ and $m_0,
\ldots, m_i$. Thus this describes a strategy for II in $G_{Y,X}$ and it
suffices to verify that it is a winning strategy.

But since II follows her strategy in $H_{Y,X}$, we know that $(w_i)\sim_C(y_i)$
and therefore, since $u_i$ and $v_i$ are defined by the same coefficients over
respectively $(y_i)$ and $(w_i)$, we have that $(v_i)\sim_C(u_i)$.
\end{proof}

\subsection{A dichotomy for minimality}
We are now in condition to prove the central result of this paper.

\begin{thm}[3rd dichotomy]\label{3rddichotomy}
Let $E$ be a Banach space with a basis $(e_i)$. Then either $E$ contains a
minimal block subspace or a continuously tight block subspace.
\end{thm}

\begin{proof}
Suppose that $E$ has no continuously tight block basic sequence. By Lemma
\ref{cont tight}, we can, modulo  passing to a block subspace, suppose that for
some constant $C$ and for all block subspaces $X\leq E$ there is a further
block subspace $Y\leq X$ such that I has no winning strategy in the game
$H_{Y,X}$ with constant $C$. By the determinacy of open games, this implies
that for all block subspaces $X\leq E$ there is a further block subspace $Y\leq
X$ such that II has a winning strategy in the game $H_{Y,X}$ with constant $C$.

A {\em state} is a pair $(a,b)$ with $a,b\in ({\bf D}'\times \F)^{<\om}$, where
$\mathbb F$ is the set of subspaces spanned by finite block sequences and ${\bf
D}'$ the set of not necessarily normalised blocks,  such that $|a|=|b|$ or
$|a|=|b|+1$. The set $S$ of states is countable, and corresponds to the
possible positions of a game $G_{L,M}$ after a finite number of moves were
made, restricted to elements that affect the outcome of the game from that
position (i.e., $m_i$'s and $n_i$'s are forgotten).

For each state $s=(a,b)$ we will define the game $G_{L,M}(s)$ in a manner
similar to the game $G_{L,M}$ depending on whether  $|a|=|b|$ or $|a|=|b|+1$.
To avoid excessive notation we do this via two examples:

If $a=(a_0,A_0,a_1,A_1)$,
$b=(b_0,B_0,b_1,B_1)$, the game $G_{L,M}(s)$ will start with II playing some integer
$n_2$, then I playing $(u_2,E_2,m_2)$ with $n_2 \leq E_2 \subseteq L$ and $u_2 \in A_0+A_1+E_2$, II playing $(v_2,F_2,n_3)$
with $m_2 \leq F_2 \subseteq M$ and $v_2 \in B_0+B_1+F_2$,  etc, and the outcome of the game will be
the pair of infinite sequences $(a_0,a_1,u_2,\ldots)$ and $(b_0,b_1,v_2,\ldots)$.

If $a=(a_0,A_0,a_1,A_1)$,
$b=(b_0,B_0)$, the game $G_{L,M}(s)$ will start with I playing some integer
$m_1$, then II playing $(v_1,F_1,n_2)$ with $m_1 \leq F_1 \subseteq M$ and $v_1 \in B_0+F_1$, I playing $(u_2,E_2,m_2)$
with $n_2 \leq E_2 \subseteq L$ and $u_2 \in A_0+A_1+E_2$,  etc, and the outcome of the game will be
the pair of infinite sequences $(a_0,a_1,u_2,\ldots)$ and $(b_0,v_1,v_2,\ldots)$.

The following lemma is well-known and easily proved by a simple diagonalisation.
\begin{lemme}
Let $N$ be a countable set and let $\mu\colon bb(E) \rightarrow \ku P(N)$ satisfy either
$$
V\leq^* W\saa \mu(V)\subseteq \mu (W)
$$
or
$$
V\leq^* W\saa \mu(V)\supseteq \mu (W).
$$
Then there exists a stabilising
block subspace $V_0 \leq E$, i.e., such that
$\mu(V)=\mu(V_0)$ for any $V \leq^* V_0$.
\end{lemme}
Let now $\tau\colon bb(E) \rightarrow \ku P(S)$ be defined by
$$
s\in \tau(M)\equi \e L\leq M \textrm{ such that player II has a
winning strategy in } G_{L,M}(s).
$$
By the asymptotic nature of the game we see that $M'\leq^* M\saa \tau(M')\subseteq \tau (M)$, and
therefore there exists $M_0\leq E$ which is stabilising for $\tau$. We
then define a map $\rho\colon bb(E)\rightarrow \ku P(S)$ by
setting
$$
s \in \rho(L)\equi \textrm { player II has a winning strategy in } G_{L,M_0}(s).
$$
Again  $L'\leq^* L\saa \rho(L')\supseteq \rho (L)$ and therefore there exists
$L_0 \leq M_0$ which is stabilising for $\rho$. Finally, the reader will easily
check that $\rho(L_0)=\tau(L_0)=\tau(M_0)$, see, e.g., \cite{anna} or
\cite{subsurfaces}.
\begin{lemme}
For every $M\leq L_0$, II has a winning strategy for the game $G_{L_0,M}$.
\end{lemme}

\begin{proof}
Fix $M$ a block subspace of $L_0$. We begin by showing that $(\tom, \tom)\in
\tau (L_0)$. To see this, we notice that as $L_0\leq E$, there is a $Y\leq L_0$
such that II has a winning strategy for $H_{Y,L_0}$ and thus, by Lemma
\ref{games}, also a winning strategy in $G_{Y,L_0}$ with constant $C$. So
$(\tom, \tom)\in \tau (L_0)$.

We will show that for all states
$$
((u_0,E_0, \ldots,u_i,E_i),(v_0,F_0,\ldots,v_i,F_i))\in \tau(L_0),
 $$
 there is an $n$ such that for all $n\leq E\subseteq L_0$ and $u\in E_0+\ldots+E_i+E$, we have
 $$
 ((u_0,E_0, \ldots,u_i,E_i,u,E),(v_0,F_0,\ldots,v_i,F_i))\in \tau(L_0).
 $$

Similarly, we show that  for all states
$$
((u_0,E_0, \ldots,u_{i+1},E_{i+1}),(v_0,F_0,\ldots,v_i,F_i))\in \tau(L_0)
$$
and for all $m$ there are $m\leq F\subseteq M$ and $v\in F_0+\ldots+F_i+F$ such that
$$
((u_0,E_0, \ldots,u_{i+1},E_{i+1}),(v_0,F_0,\ldots,v_i,F_i,v,F))\in \tau(L_0).
$$
Since the winning condition of $G_{L_0,M}$ is closed, this clearly shows that
II has a winning strategy in $G_{L_0,M}$ (except for the integers $m$ and $n$,
$\tau(L_0)$ is a winning quasi strategy for II).

So suppose that
$$
s=((u_0,E_0, \ldots,u_i,E_i),(v_0,F_0,\ldots,v_i,F_i))\in \tau(L_0)=\rho(L_0),
$$
then II has a winning strategy in $G_{L_0,M_0}(s)$ and hence there is an $n$
such that for all $n\leq E\subseteq L_0$ and $u\in E_0+\ldots+E_i+E$, II has a
winning strategy in $G_{L_0,M_0}(s')$, where
$$
s'= ((u_0,E_0, \ldots,u_i,E_i,u,E),(v_0,F_0,\ldots,v_i,F_i)).
$$
So $s'\in \rho(L_0)=\tau(L_0)$.

Similarly, if
$$
s=((u_0,E_0, \ldots,u_{i+1},E_{i+1}),(v_0,F_0,\ldots,v_i,F_i))\in \tau(L_0)=\tau(M)
$$
and $m$ is given, then as II has a winning  strategy for $G_{L,M}(s)$ for some
$L\leq M$, there are  $m\leq F\subseteq M$ and $v\in F_0+\ldots+F_i+F$ such
that II has a winning strategy in $G_{L,M}(s')$, where
$$
s'=((u_0,E_0, \ldots,u_{i+1},E_{i+1}),(v_0,F_0,\ldots,v_i,F_i,v,F)).
$$
So $s'\in \tau(M)=\tau(L_0)$.
\end{proof}
Choose now $Y=[y_i]\leq L_0$ such that II has a winning strategy in
$H_{Y,L_0}$. We shall show that any block subspace $M$ of $L_0$ contains a
$C^2$-isomorphic copy of $Y$, which implies that $Y$ is $C^2+\epsilon$-minimal
for any $\epsilon>0$.

To see this, notice that, since II has a winning strategy in $H_{Y,L_0}$,
player I has a strategy in the game $G_{L_0,M}$ to produce a sequence $(u_i)$
that is $C$-equivalent with the basis $(y_i)$. Moreover, we can ask that I
plays $m_i=0$. Using her winning strategy for $G_{L_0,M}$, II can then respond
by producing a sequence $(v_i)$ in $M$ such that $(v_i)\sim_C(u_i)$. So
$(v_i)\sim_{C^2}(y_i)$ and $Y\sqsubseteq_{C^2} M$.
\end{proof}

\

Finally we observe that by modifying the notion of embedding in the definition
of a tight basis, we obtain variations of our dichotomy theorem with a weaker
form of tightness on one side and a stronger form of minimality on the other.

\begin{thm}
Every Banach space with a basis contains a block subspace $E=[e_n]$ which satisfies one of the two following properties:
\begin{enumerate}
\item For any $[y_i] \leq E$, there exists a sequence $(I_i)$ of successive
intervals such that for any infinite subset $A$ of $\N$, the basis $(y_i)$
does not embed into $[e_n]_{n \notin \cup_{i \in A}I_i}$ as a sequence of
disjointly supported blocks , resp. as a permutation of a block sequence,
resp. as a block sequence.
\item For any $[y_i] \leq E$,  $(e_n)$ is equivalent to a sequence of disjointly
supported blocks of $[y_i]$, resp. $(e_n)$ is permutatively equivalent to a
block sequence of $[y_i]$, resp. $(e_n)$ is equivalent to a block sequence
of $[y_i]$.
\end{enumerate}
\end{thm}
The case of block sequences immediately implies the theorem of Pe\l czar \cite{anna}.

The fact that the canonical basis of $T^*$ is strongly asymptotically
$\ell_{\infty}$  implies easily that it is tight for ``embedding as a sequence
of disjointly supported blocks'' although $T^*$ is minimal in the usual sense.
We do not know of other examples of spaces combining one form of  minimality
with another form of  tightness in the above list.

\section{Tightness with constants and crude stabilisation of local structure}
We shall now consider a stronger notion of tightness, which is essentially
local in nature. Let $E$ be a space with a basis $(e_n)$. There is a
particularly simple case when a sequence $(I_i)$ of intervals associated to a
subspace $Y$ characterises the tightness of $Y$ in $(e_n)$. This is when for
all integer constants $K$, $Y \not\sqsubseteq_K [e_n]_{ n \notin I_K}$. This
property has the following useful reformulations.
\begin{prop}\label{tightnesswithconstants} Let $E$ be a space with a basis $(e_n)$. The following are equivalent:
 \begin{enumerate}
\item For any block sequence $(y_n)$  there are intervals $I_0<I_1<I_2<\ldots$ such that for all $K$,
$$
[y_n]_{n\in I_K}\not\sqsubseteq_K[e_n]_{n\notin I_K}.
$$
\item  For any space $Y$, there are intervals $I_0<I_1<I_2<\cdots$ such that for all $K$,
$$
Y \not\sqsubseteq_K[e_n]_{ n \notin I_K}.
$$
\item  No space embeds uniformly into the tail subspaces of $E$.
\item  There is no $K$ and no subspace of $E$ which is $K$-crudely finitely representable in any tail subspace of $E$.
\end{enumerate}
\end{prop}

A basis satisfying properties (1), (2), (3), (4), as well as the space it
generates, will be said to be {\em tight with constants}.

\begin{proof} The implications (1)$\saa$(2)$\saa$(3)  are clear.

To prove (3)$\saa$(4) assume some subspace $Y$ of $E$ is $K$-crudely finitely
representable in any tail subspace of $E$. Without loss of generality, we may
assume that $Y=[y_n]$ is a block subspace of $E$. We pick a subsequence $(z_n)$
of $(y_n)$ in the following manner. Let $z_0=y_0$, and if $z_0,\ldots,z_{k-1}$
have been chosen, we choose $z_k$ supported far enough on the basis $(e_n)$, so
that $[z_0,\ldots,z_{k-1}]$ has a $2K$-isomorphic copy in  $[e_n \del k \leq n
<\min({\rm supp }\; z_k)]$. It follows that for any $k$, $Z=[z_n]$ has an
$M$-isomorphic copy in the tail subspace $[e_n \del n \geq k]$ for some $M$
depending only on $K$ and the constant of the basis $(e_n)$.

To prove (4)$\saa$(1), let $c(L)$ be a constant such that if two block
sequences differ in at most $L$ terms, then they are $c(L)$-equivalent. Now
assume (4) holds and let $(y_n)$ be a block sequence of $(e_n)$. Suppose also
that $I_0<\ldots<I_{K-1}$ have been  chosen. By (4) applied to $Y=[y_n]_{n=\max
I_{K-1}+1}^\infty$, we can then find $m$ and $l>\max I_{K-1}$ such that
$[y_n]_{n=\max I_{K-1}+1}^l$ does not $K\cdot c(\max I_{K-1}+1)$-embed into
$[e_n]_{n=m}^\infty$.  Let now
$$
I_K=[\max I_{K-1}+1,l+m]
$$
and notice that, as $[y_n]_{n=\max I_{K-1}+1}^l\subseteq [y_n]_{n\in I_K}$, we
have that $ [y_n]_{n\in I_K}$ does not $K\cdot c(\max I_{K-1}+1)$-embed into
$[e_n]_{n=m}^\infty$. Also, since $(e_n)_{n=m}^\infty$ and
$$
(e_n)_{n=0}^{\max I_{K-1}}{}^\frown(e_n)_{n=\max I_{K-1}+1+m}^\infty
$$
only differ in $\max I_{K-1}+1$ many terms, $[y_n]_{n\in I_K}$ does not $K$-embed into
$$
[e_n]_{n=0}^{\max I_{K-1}}+[e_n]_{n=\max I_{K-1}+1+m}^\infty,
$$
and thus not into the subspace $[e_n]_{n\notin I_K}$ either.
\end{proof}
It is worth noticing that a basis $(e_n)$, tight with constants, is necessarily
continuously tight. For a simple argument shows that in order to find the
intervals $I_K$ satisfying (1) above, one only needs to know a beginning of the
block sequence $(y_n)$ and hence the intervals can be found continuously in
$(y_n)$. From Proposition \ref{tightnesswithconstants} we also deduce that any
block basis or shrinking basic sequence in the span of a tight with constants
basis is again tight with constants.

There is a huge difference between the fact that no {\em subspace} of $E$ is
$K$-crudely finitely representable in all tails of $E$ and then that no {\em space}
is $K$-crudely finitely representable in all tails of $E$. For example, we
shall see that while the former holds for Tsirelson's space, by Dvoretzky's
Theorem (see e.g. \cite{FLM}), $\ell_2$ is always finitely representable in any Banach space.

Recall that a basis $(e_n)$  is said to be {\em strongly asymptotically
$\ell_p$}, $1 \leq p \leq +\infty$, \cite{DFKO}, if there exists a constant $C$
and a function $f:\N \rightarrow \N$ such that for any $n$, any family of $n$
unit vectors which are disjointly supported in $[e_k \del k \geq f(n)]$ is
$C$-equivalent to the canonical basis of $\ell_p^n$.

\begin{prop}\label{dfko}
Let $E$ be a Banach space with a strongly asymptotically $\ell_p$ basis
$(e_n)$, $1 \leq p<+\infty$, and not containing a copy of $\ell_p$. Then
$(e_n)$ is tight with constants.\end{prop}

\begin{proof}
Assume that some Banach space $Y$ embeds with constant $K$ in any tail subspace
of $E$. We may assume that $Y$ is generated by a block-sequence $(y_n)$ of $E$
and, since any strongly asymptotically $\ell_p$ basis is unconditional, $(y_n)$
is unconditional. By renorming $E$ we may assume it is $1$-unconditional. By a
result of W.B. Johnson \cite{J} for any $n$ there is a constant $d(n)$ such
that $(y_0,\ldots,y_n)$ is $2K$-equivalent to a sequence of vectors in the
linear span of $d(n)$ disjointly supported unit vectors in any tail subspace of
$E$, in particular in $[e_k \del k \geq f(d(n))]$. Therefore $[y_0,\ldots,y_n]$
$2KC$-embeds into $\ell_p$. This means that $Y$ is crudely finitely
representable in $\ell_p$ and therefore embeds into $L_p$, and since $(y_n)$ is
unconditional asymptotically $\ell_p$, that $Y$ contains a copy of $\ell_p$
(details of the last part of this proof may be found in \cite{DFKO}).
\end{proof}

\begin{cor}\label{tsi}
Tsirelson's space $T$ and its convexifications $T^{(p)}$, $1<p<+\infty$, are
tight with constants.
\end{cor}

Observe that on the contrary, the dual $T^*$ of $T$, which is strongly
asymptotically $\ell_{\infty}$ and does not contain a copy of $c_0$, is minimal
and therefore does not contain any tight subspace.

\

Suppose a space $X$ is crudely finitely representable in all of its subspaces.
Then there is some constant $K$ and a subspace $Y$ such that $X$  is
$K$-crudely finitely representable in all of the subspaces of $Y$. For if not,
we would be able to construct a sequence of basic sequences $(x_n^K)$ in $X$
such that $(x_n^{K+1})$ is a block sequence of $(x_n^K)$ and such that $X$ is
not $K^2$-crudely finitely representable in $[x_n^{K}]$. By Lemma \ref{block
uniformisation}, we can then find a block sequence $(y_n)$ of $(x_n^0)$ that is
$\sqrt K$-equivalent with a block sequence of $(x_n^K)$ for any $K$ and hence
if $X$ were $K$-crudely finitely representable in $[y_n]$ for some $K$, then it
would also be $K^{3/2}$-crudely finitely representable in $[x_n^{K}]$, which is
a contradiction.

When a space $X$ is $K$-crudely finitely representable in any of its subspaces
for some $K$, we say that $X$ is {\em locally minimal}. For example, by the
universality properties of $c_0$, any space
with an asymptotically $\ell_{\infty}$ basis is locally minimal.

\begin{thm}[5th dichotomy]\label{5th}
Let $E$ be an infinite-dimensional Banach space with basis $(e_n)$. Then there
is a block sequence $(x_n)$ satisfying one of the following two properties,
which are mutually exclusive and both possible.
\begin{enumerate}
 \item $(x_n)$ is tight with constants,
  \item $[x_n]$ is locally minimal.
\end{enumerate}
\end{thm}

\begin{proof}
If $E$ contains $c_0$, the result is trivial. So suppose not and find by the
solution to the distortion problem a block sequence $(y_n)$ and inevitable,
positively separated, closed subsets $F_0$ and $F_1$ of the unit sphere of
$[y_n]$. Define for each integer $K\geq 1$ the set
\[\begin{split}
\A_K=\{& (z_n)\leq (y_n)\del z_{2n}\in F_0\cup F_1 \textrm{ and $(z_{2n})$ codes by $0$'s and $1$'s a block}\\
 &\textrm{sequence $(v_n)$ of $(z_{2n+1})$ such that for all }N, [v_n]\sqsubseteq_K [z_{2n+1}]_{n\geq N}\\
&\textrm { and moreover } 1/2<\|v_n\|<2\}.
\end{split}\]
Clearly $\A_K$ is analytic, so we can apply Gowers' Determinacy Theorem to get one of two cases
\begin{itemize}
  \item [(i)] either there is a block sequence $(x_n)$ and a $K$ such that player II has
  a strategy to play inside $(\A_K)_\Delta$ whenever I plays a block
  sequence of $(x_n)$, where $\Delta$ will be determined later,
  \item [(ii)] or we can choose inductively a sequence of block sequences $(x_n^K)$ such
  that $(x_n^{K+1})\leq (x_n^K)$ and such that no block sequence of
  $(x_n^K)$ belongs to $\A_K$.
\end{itemize}

Consider first case (ii). Set $w_n=x^n_{2n}$ and choose now  further block
sequences $(x_n)$ and $(h_n)$ of $(w_n)$ such that
$$
x_0<h_0<h_1<x_1<h_2<h_3<x_4<\ldots
$$
and $h_{2n}\in F_0$, $h_{2n+1}\in F_1$.

We claim that $(x_n)$ is tight with constants. If not, we can find some block
sequence $(u_n)$ of $(x_n)$  and a $K$ such that $[u_n]$ embeds with constant
$K$ into any tail subspace of $[x_n]$. By passing to tails of $(x_n)$ and of
$(u_n)$, we can suppose that $(x_n)$ is a block sequence of $(x_n^{K})$,
$(u_n)$ is a block sequence of $(x_n)$ and $[u_n]$ $K$-embeds into all tails of
$[x_n]$. By filling in with appropriate $h_i$ between $x_n$ and $x_{n+1}$, we
can now produce a block sequence $(z_n)$ of $(x_n^{K})$ such that $(z_{2n})$
codes by $0$'s and $1$'s the block sequence $(u_n)$ of $(z_{2n+1})$ with the
property that for all $N$, $[u_n]\sqsubseteq_K [z_{2n+1}]_{n\geq N}$. In other
words, we have produced a block sequence of $(x_n^K)$ belonging to $\A_K$,
which is impossible. Thus, $(x_n)$ is tight with constants.

Consider now case (i) instead and let II play according to her strategy. We
suppose that $\Delta$ is chosen sufficiently small so that $\delta_i<{\rm
dist}(F_0,F_1)/3$ and if two block sequences are $\Delta$-close then they are
$2$-equivalent. Let $(y_n)\in (\A_K)_\Delta$ be the response by II to the
sequence $(x_n)$ played by I and let $(z_n)\in \A_K$ be such that
$\|z_n-y_n\|<\delta_n$ for all $n$. Then $(z_{2n})$ codes by $0$'s and $1$'s a
block sequence $(v_n)$ of $(z_{2n+1})$. Let $(u_n)$ be the block sequence of
$(y_{2n+1})$ constructed in the same way as $(v_n)$ is constructed over
$(z_{2n+1})$. We claim that $(u_n)$ is $4K$-crudely finitely representable in
any block subspace of $[x_n]$.

For this, let $[u_0,\ldots,u_m]$ be given and suppose that $(f_n)$ is any block
subspace of $(x_n)$. Find a large $k$ such that $(z_0,z_2,\ldots,z_{2k})$ codes
the block sequence $(v_0,\ldots, v_m)$ of $(z_1,\ldots,z_{2k+1})$ and let $l$
be large enough so that when I has played $x_0,\ldots, x_l$ then II has played
$y_0,\ldots, y_{2k+1}$. Consider now the game in which player I plays
$$
x_0,x_1,\ldots,x_l,f_{l+1},f_{l+2},\ldots.
$$
Then, following the strategy, II will play a block sequence
$$
y_0,\ldots, y_{2k+1},g_{2k+2},g_{2k+3},\ldots\in (\A_K)_\Delta.
$$
So let $(h_n)\in \A_K$ be such that  $\|h_n-y_n\|<\delta_n$ for all $n\leq
2k+1$ and $\|h_n-g_n\|<\delta_n$ for all $n\geq 2k+2$. For $n\leq k$, we have,
as $\|h_{2n}-z_{2n}\|<2\delta_n<\frac23{\rm dist}(F_0,F_1)$, that $h_{2n}\in
F_i\equi z_{2n}\in F_i$.  Also, $(h_{2n+1})_{n=0}^k$ and $(y_{2n+1})_{n=0}^k$
are $2$-equivalent and $(h_{2n+1})_{n=k+1}^\infty$ and
$(g_{2n+1})_{n=k+1}^\infty$ are $2$-equivalent, so $(h_{2n})$ will code a block
sequence $(w_n)$ of $(h_{2n+1})$ such that $(w_0,\ldots,w_m)$ is $2$-equivalent
to $(u_0,\ldots, u_m)$. Moreover, since $(h_n)\in \A_K$, $[w_n]$ will $K$-embed
into every tail subspace of $[h_{2n+1}]$, and hence $2K$-embed into every tail
subspace of $[g_{2n+1}]$. Therefore, since $(g_{2n+1})$ is block sequence of
$(f_n)$, $[u_0,\ldots,u_m]$ will $4K$-embed into $[f_n]$, which proves the
claim. It follows that $[u_n]$ is locally minimal, which proves the theorem.
\end{proof}

Local minimality can be reformulated in a way that makes the relation to local
theory clearer. For this, let $\F_n$ be the metric space of all $n$-dimensional
Banach spaces up to isometry equipped with the Banach-Mazur metric
$$
d(X,Y)=\inf\big(\log( \|T\|\cdot\|T\inv\|)\del T\colon X\til Y \textrm{ is an isomorphism }\big).
$$
Then for every Banach space $X$, the set of $n$-dimensional $Y$ that are almost
isometrically embeddable into $X$ form a closed subset $(X)_n$ of $\F_n$. It is
well-known that this set $(X)_n$ does not always stabilise, i.e., there is not
necessarily a subspace $Y\subseteq X$ such that for all further subspaces
$Z\subseteq Y$, $(Z)_n=(Y)_n$. However, if instead $X$ comes equipped with a
basis and for all block subspaces $Y$ we let $\{Y\}_n$ be the set of all
$n$-dimensional spaces that are almost isometrically embeddable into all tail
subspaces of $Y$, then one can easily stabilise $\{Y\}_n$ on subspaces. Such
considerations are for example the basis for \cite{MMT}.

Theorem \ref{5th} gives a dichotomy for when one can stabilise the set $(X)_n$
in a certain way, which we could call {\em crude}. Namely, $X$ is locally
minimal if and only if there is some constant $K$ such that for all subspaces
$Y$ of $X$ and all $n$,  $d_H\big((X)_n,(Y)_n\big)\leq K$, where $d_H$ is the
Hausdorff distance. So by Theorem \ref{5th}, the local structure stabilises
crudely on a subspace if and only if a space is not saturated by basic
sequences tight with constants.

\

Often it is useful to have a bit more than local minimality. So we say that a
basis $(e_n)$ is {\em locally block minimal} if it is $K$-crudely block
finitely representable in all of its block bases for some $K$. As with crude
finite representability we see that there then must be a constant $K$ and a
block $(y_n)$ such that $(e_n)$ is $K$-crudely block finitely representable in
all  block subspaces of $(y_n)$. We now have the following version of Theorem
\ref{5th} for finite block representability.

\begin{thm}\label{5th block}
Let $(e_n)$ be a Schauder basis. Then $(e_n)$ has a block basis $(x_n)$ with
one of the following two properties, which are mutually exclusive and both
possible.
\begin{enumerate}
 \item For all block bases $(y_n)$ of $(x_n)$ there are intervals
 $I_1<I_2<I_3<\ldots$ such that $(y_n)_{ n\in I_K}$ is not $K$-equivalent
 to a block sequence of $(x_n)_{ n\notin I_K}$,
  \item $(x_n)$ is locally block minimal.
\end{enumerate}
\end{thm}

\

Finally we note that there exist tight spaces which do not admit subspaces
which are tight with constants:

\begin{ex}\label{OdellSchlumprecht}
There exists a reflexive, tight, locally block minimal Banach space.
\end{ex}

\begin{proof}
E. Odell and T. Schlumprecht \cite{OS:universalhi}  have built a reflexive
space $OS$ with a basis such that every monotone basis is block finitely
representable in any block subspace of $OS$. It is in particular locally block
minimal and therefore contains no basic sequence which is tight with constants.
We do not know whether the space $OS$  is tight. Instead, we notice that since
the summing basis of $c_0$ is block finitely representable in any block
subspace of $OS$, $OS$ cannot contain an unconditional block sequence. By
Gowers' 1st dichotomy it follows that some block subspace of $OS$ is HI, and by
the 3rd dichotomy (Theorem \ref{3rddichotomy}) and the fact that HI spaces do
not contain minimal subspaces, that some further block subspace is tight, which
completes the proof.
\end{proof}

It is unknown whether there is an unconditional example with the above
property. There exists an unconditional version of $OS$ \cite{OS:universalunc},
but it is unclear whether it has no minimal subspaces. However, the dual of a
space constructed by Gowers in \cite{g:hyperplanes} can be shown to be both
tight and locally minimal.
\begin{ex}\label{gunc}\cite{exemples}
There exists a space with an unconditional basis which is tight and locally
minimal.
\end{ex}

\section{Local block minimality, asymptotic structure and a dichotomy for containing $c_0$ or $\ell_p$}
Recall that a basis $(e_n)$ is said to be {\em asymptotically $\ell_p$} (in the
sense of Tsirelson's space) if there is a constant $C$
such that for all normalised block sequences $n<x_1<\ldots<x_n$,
$(x_i)_{i=1}^n$ is $C$-equivalent with the standard unit vector basis of
$\ell_p^n$.

When $(x_n)$ is asymptotically $\ell_p$ and some block subspace $[y_n]$ of
$[x_n]$ is $K$-crudely  block finitely representable in all tail subsequences
of $(x_n)$, then it is clear that  $(y_n)$ must actually be equivalent with the
unit vector basis of $\ell_p$, or $c_0$ for $p=\infty$. So this shows that for
asymptotically $\ell_p$ bases $(x_n)$, either $[x_n]$ contains an isomorphic
copy of $\ell_p$ or $c_0$ or $(x_n)$ itself satisfies  condition (1) of Theorem
\ref{5th block}. This is the counterpart of Proposition \ref{dfko} for block
sequences. As an example, we mention that, since $T^*$ does not contain $c_0$,
but has a strongly asymptotically $\ell_\infty$ basis, it thus satisfies (1).
This small observation indicates that one can characterise when $\ell_p$ or
$c_0$ embeds into a Banach space by characterising when a space contains an
asymptotic $\ell_p$ space.

We first prove a dichotomy for having an asymptotic $\ell_p$ subspace, for the
proof of which we need the following lemma.

\begin{lemme}\label{krivine}
Suppose $(e_n)$ is a basic sequence such that for some $C$ and all $n$ and
normalised block sequences  $n<y_1<y_2<\ldots<y_{2n}$, we have
$$
(y_{2i-1})_{i=1}^n\sim_C(y_{2i})_{i=1}^n.
$$
Then $(e_n)$ has an asymptotic $\ell_p$ subsequence for some $1\leq p\leq
\infty$.
\end{lemme}

\begin{proof}
By the Theorem of Brunel and Sucheston \cite{brunel}, we can, by passing to a
subsequence of $(e_n)$, suppose that $(e_n)$ generates a spreading model, i.e.,
we can assume that for all integers $n<l_1<l_2<\ldots<l_n$ and
$n<k_1<k_2<\ldots<k_n$ we have
$$
(e_{l_1},\ldots,e_{l_n})\sim_{1+\frac 1n}(e_{k_1},\ldots,e_{k_n}).
$$
Now suppose that $e_n<y_1<y_2<\ldots<y_{n}$ and $e_n<z_1<\ldots<z_n$ are
normalised block sequences of $(e_{2i})$. Then there are $n<l_1<l_2<\ldots<l_n$
and $n<k_1<k_2<\ldots<k_n$ such that
$$
e_n<y_1<e_{l_1}<y_2<e_{l_2}<\ldots<y_n<e_{l_n}
$$
and
$$
e_n<z_1<e_{k_1}<z_2<e_{k_2}<\ldots<z_n<e_{k_n},
$$
so $(y_i)\sim_C(e_{l_i})\sim_{1+\frac1n}(e_{k_i})\sim_C(z_i)$ and
$(y_i)\sim_{2C^2}(z_i)$. Thus, asymptotically all finite normalised block
sequences are $2C^2$-equivalent.

By Krivine's Theorem \cite{krivine}, there is some $\ell_p$ that is block
finitely representable in $(e_{2i})$ and hence asymptotically all finite
normalised block sequences are equivalent to $\ell_p^n$ of the correct
dimension and $(e_{2i})$ is asymptotic $\ell_p$.
\end{proof}

\begin{thm}\label{asymp}
Suppose $X$ is a Banach space with a basis. Then $X$ has a block subspace $W$,
which is either asymptotic $\ell_p$, for some $1\leq p\leq +\infty$, or such
that
$$
\a M \;\e n\;  \a U_1,\ldots, U_{2n}\subseteq W\; \e u_i\in \ku S_{U_i}\;\Big(u_1<\ldots<u_{2n}\;\&\:
 (u_{2i-1})_{i=1}^n\not \sim_M(u_{2i})_{i=1}^n\Big).
$$
\end{thm}

\begin{proof}
Assume first that for some $M$ and $V\subseteq X$ we have
$$
\a n\;\a Y\subseteq V\;  \e Z\subseteq Y\; \a z_1<\ldots<z_{2n}\in \ku S_Z\; (z_{2i-1})_{i=1}^n\sim_M(z_{2i})_{i=1}^n.
$$
Then we can inductively define $V\supseteq Z_1\supseteq Z_2\supseteq
Z_3\supseteq\ldots$ such that for each $n$,
$$
\a z_1<\ldots<z_{2n}\in \ku S_{Z_n}\; (z_{2i-1})_{i=1}^n\sim_M(z_{2i})_{i=1}^n.
$$
Diagonalising over this sequence, we can find a block subspace $W=[w_n]$ such
that for all $m\geq n$, $w_m\in Z_n$. Therefore, if $n<z_1<\ldots<z_{2n}$ is a
sequence of normalised blocks of $(w_n)$, then
$(z_{2i-1})_{i=1}^n\sim_M(z_{2i})_{i=1}^n$. By Lemma \ref{krivine}, it follows
that $W$ has an asymptotic $\ell_p$ subspace.

So suppose on the contrary that
$$
\a M\; \a V\subseteq X\; \e n\;\e Y\subseteq V\; \a Z\subseteq Y\;
\e z_1<\ldots<z_{2n}\in \ku S_Z\; (z_{2i-1})_{i=1}^n\not\sim_M(z_{2i})_{i=1}^n.
$$
Now find some small $\Delta>0$ such that two normalised block sequences of $X$
that are $\Delta$-close are $\sqrt2$-equivalent. Applying the determinacy
theorem of Gowers to $\Delta$ and the sets
$$
\A(M,n)=\{(z_i)\del (z_{2i-1})_{i=1}^n\not\sim_M(z_{2i})_{i=1}^n\},
$$
we have that
\begin{displaymath}\begin{split}
\a M\; \a V\subseteq X\;\e n\; \e Z\subseteq V\quad &\textrm{II has a strategy
to play normalised}\\
&\textrm{$z_1<\ldots<z_{2n}$ such that
$(z_{2i-1})_{i=1}^n\not\sim_\frac M2(z_{2i})_{i=1}^n$}\\
&\textrm{when I is restricted to playing blocks of $Z$.}
\end{split}\end{displaymath}

Using this we can inductively define a sequence $X\supseteq W_1\supseteq
W_2\supseteq W_3\supseteq\ldots$ such that for all $M$ there is an $n=n(2M)$
such that II has a strategy to play normalised $z_1<\ldots<z_{2n}$ such that
$(z_{2i-1})_{i=1}^n\not\sim_M(z_{2i})_{i=1}^n$ whenever  I is restricted to
playing blocks of $W_M$. Diagonalising over this sequence, we find some
$W\subseteq X$ such that for all $M$, $W\subseteq^*W_M$. So for all $M$ there
is $n$ such that II has a strategy to play normalised $z_1<\ldots<z_{2n}$ such
that $(z_{2i-1})_{i=1}^n\not\sim_M(z_{2i})_{i=1}^n$ whenever  I is restricted
to playing blocks of $W$.

Letting player I play a segment of the block basis of $U_i$ until II plays a
vector, we easily see that whenever $U_1,\ldots,U_{2n}\subseteq W$, there are
$z_i\in \ku S_{U_i}$ such that $z_1<\ldots<z_{2n}$ and
$(z_{2i-1})_{i=1}^n\not\sim_M(z_{2i})_{i=1}^n$. This finishes the proof.
\end{proof}

Using this, we can now prove the main result of this section.

\begin{thm}[The $c_0$ and $\ell_p$ dichotomy]
Suppose $X$ is a Banach space not containing a copy of $c_0$ nor of $\ell_p$,
$1\leq p<\infty$. Then $X$ has a subspace $Y$ with a basis satisfying one of
the following properties.
\begin{itemize}
  \item [(i)] $\a M \;\e n\;  \a U_1,\ldots, U_{2n}\subseteq Y\; \e u_i\in \ku S_{U_i}\;u_1<\ldots<u_{2n}\;\&\:
 (u_{2i-1})_{i=1}^n\not\sim_{M}(u_{2i})_{i=1}^n$.
  \item [(ii)] For all block bases $(z_n)$ of $Y=[y_n]$ there are
intervals $I_1 < I_2 < I_3 <\ldots$ such that $(z_n)_{n\in I_K}$ is not
$K$-equivalent to a block sequence of $(y_n)_{n\notin I_K}$.
\end{itemize}
\end{thm}
Notice that both (i) and (ii) trivially imply that $Y$ cannot contain a copy of
$c_0$ or $\ell_p$, since (i) implies some lack of homogeneity and (ii) some
lack of minimality. However, we do not know if there are any spaces satisfying
both (i) and (ii) or if, on the contrary, these two properties are
incompatible.

\begin{proof}If $X$ has no block subspace satisfying (i), then it must have
an asymptotically $\ell_p$ block subspace $Y$ for some $1\leq p\leq \infty$.
Let $C$ be the constant of asymptoticity. Suppose now that $Z=[z_n]$ is a
further block subspace. Since $Z$ is not isomorphic to $\ell_p$, this means
that for all $L$ and $K$, there is some $M$ such that $(z_n)_{n=L}^M$ is not
$KC$-equivalent to $\ell_p^{M-L-1}$ and hence not $K$-equivalent with a
normalised block sequence of $(y_n)_{n=M}^\infty$ either. So if
$I_1<I_2<\ldots<I_{K-1}$ have been defined, to define $I_K$, we let $N=\max
I_{K-1}+1$ and find $M$ such that $(z_n)_{n=2N-1}^M$ is not $K$-equivalent with
a normalised block sequence of $(y_n)_{n=M}^\infty$. It follows that
$(z_n)_{n=N}^M$ is not $K$-equivalent with a normalised block sequence of
$(y_1,y_2,\ldots,y_{N-1},y_{M+1},y_{M+2},\ldots)$. Letting $I_K=[N,M]$ we have
the result.
\end{proof}

We should mention that G. Androulakis, N. Kalton and Tcaciuc \cite{T2} have extended Tcaciuc's dichotomy from \cite{T} to a
dichotomy characterising containment of $\ell_p$ and $c_0$. The result above implies theirs and moreover provides additional
information.

\section{Tightness by range and subsequential minimality}
Theorem \ref{main} shows that if one allows oneself to pass to a basis for a
subspace, one can find a basis in which there is  a close connection between
subspaces spanned by block bases and subspaces spanned by subsequences. Thus,
for example, if the basis is tight there can be no space embedding into all the
subspaces spanned by subsequences of the basis. On the other hand, any block
basis in Tsirelson's space $T$ is equivalent to a subsequence of the basis, and
actually every subspace of a block subspace $[x_n]$ in $T$  contains an
isomorphic copy of a subsequence of $(x_n)$.  In fact, this phenomenon has a
deeper explanation and we shall now proceed to show that the connection between
block sequences and subsequences  can be made even closer.

\begin{lemme}\label{flat subspaces}
If $(e_n)$ is a basis for a space not containing  $c_0$, then for all finite
intervals $(I_n)$ such that $\min I_n\Lim{n\til \infty} \infty$ and all
subspaces $Y$, there is a further subspace $Z$ such that
$$
\norm{P_{I_k}|_ Z}\Lim{k\til \infty}0.
$$
\end{lemme}

\begin{proof}
By a standard perturbation argument, we can suppose that $Y$ is generated by a
normalised block basis $(y_n)$. Let $K$ be the basis constant of $(e_n)$. As
$\min I_n\Lim{n\til \infty}\infty$ and each $I_n$ is finite, we can choose a
subsequence $(v_n)$ of $(y_n)$ such that for all $k$ the interval $I_k$
intersects the range of at most one vector $v_m$ from $(v_n)$. Now, since $c_0$
does not embed into $[e_n]$, no tail sequence of $(v_n)$ can satisfy an upper
$c_0$ estimate. This implies that for all $N$ and $\delta>0$ there is a
normalised vector
$$
z=\sum_{i=N}^{N'}\eta_iv_i,
$$
where $|\eta_i|<\delta$. Using this, we now construct a normalised block
sequence $(z_n)$ of $(v_n)$ such that  there are $m(0)<m(1)<\ldots$ and
$\alpha_i$ with
$$
z_n=\sum_{i=m(n)}^{m(n+1)-1}\alpha_iv_i
$$
and $|\alpha_i|<\frac 1n$ whenever $m(n)\leq i <m(n+1)$.

Now suppose $u=\sum_j\lambda_jz_j$ and $k$ are given. Then there is at most one
vector $z_n$ whose range intersect the interval $I_k$. Also, there is  at most
one vector $v_p$ from the support of $z_n$ whose range intersect $I_k$.
Therefore,
\begin{displaymath}\begin{split}
\|P_{I_k}(u)\|&=\|P_{I_k}(\lambda_nz_n)\|= \|P_{I_k}(\lambda_n\alpha_pv_p)\|\\
&\leq 2K\|\lambda_n\alpha_pv_p\|\leq |\lambda_n|\cdot\frac {2K}n\leq \frac {4K^2}n \|u\|.
\end{split}\end{displaymath}
It follows that $\|P_{I_k}|_{[z_l]}\|\leq \frac{4K^2}{n_k}$, where $n_k$ is
such that $I_k$ intersects the range of $z_{n_k}$ (or $n_k=k$ if $I_k$
intersects the range of no $z_n$). Since $\min I_k\Lim{k\til \infty}\infty$ and
$(z_n)$ is a block basis, $n_k\til \infty$ when $k\til \infty$, and hence
$\|P_{I_k}|_{[z_l]}\|\Lim{k\til \infty} 0$.
\end{proof}
Our next result should be contrasted with the construction by Pe\l czy\'nski
\cite{pe} of a basis $(f_i)$ such that every basis is equivalent with a
subsequence of it, and hence such that every space contains an isomorphic copy
of a subsequence. We shall see that for certain spaces $E$ such constructions
cannot be done relative to the subspaces of $E$ provided that we demand that
$(f_n)$ lies in $E$ too. Recall that two Banach spaces are said to be {\em
incomparable} if  neither of them embeds into the other.

\begin{prop}\label{opposite}
Suppose that $(e_n)$ is a basis such that any two block subspaces with disjoint
ranges are incomparable. Suppose also that $(f_n)$ is either a block basis or a
shrinking basic sequence in $[e_n]$. Then $[e_n]$ is saturated with subspaces
$Z$ such that no subsequence of $(f_n)$ embeds into $Z$.
\end{prop}

\begin{proof}Let $Y$ be an arbitrary subspace of $[e_n]$.
Suppose first that $(f_n)$ is a normalised shrinking basic sequence. Then, by
taking a perturbation of $(f_n)$, we can suppose that each $f_n$ is a finite
block vector of $(e_i)$ and, moreover, that $\min{\rm range}(f_n)\til \infty$.
Let $I_n={\rm range}(f_n)$.

Fix an infinite set $N\subseteq \N$. Then for all infinite subsets $A\subseteq
N$ there is an infinite subset $B\subseteq A$ such that $(f_n)_{n\in B}$ is a
block sequence and hence, since  block subspaces of $(e_n)$ with disjoint
ranges are incomparable,  $[f_n]_{n\in B}\not\sqsubseteq [e_n]_{n\notin
\bigcup_{i\in B}I_i}$, and so also $[f_n]_{n\in N}\not\sqsubseteq
[e_n]_{n\notin \bigcup_{i\in A}I_i}$. Applying Lemma \ref{projection tightness}
to $X=[f_n]_{n\in N}$ and $(I_n)_{n\in N}$, this implies that for all
embeddings $T\colon [f_n]_{n\in N}\til [e_n]_{n\in \N}$, we have $\liminf_{n\in
N}\norm{P_{I_k}T}>0$. So find by Lemma \ref{flat subspaces} a subspace
$Z\subseteq Y$ such that $\norm{P_{I_k}|_ Z}\Lim{k\til \infty}0$. Then no
subsequence of $(f_n)_{\in \N}$ embeds into $Z$.

The argument in the case $(f_n)$ is a block basis is similar.
We set $I_n={\rm range}\;f_n$ and repeat the argument above.
\end{proof}

We notice that in the above proof we actually have a measure for how ``flat'' a
subspace $Z$ of  $[e_n]$ needs to be in order that the subsequences of $(f_n)$
cannot embed into $Z$. Namely, it suffices that $\norm{P_{I_k}|_ Z}\Lim{k\til
\infty}0$.

We should also mention that, by similar but simpler arguments, one can show
that if $(e_n)$ is a basis such that any two disjoint subsequences span
incomparable spaces, then some subspace of $[e_n]$ fails to contain any
isomorphic copy of a subsequence of $(e_n)$.

The assumption in Proposition \ref{opposite} that block subspaces with disjoint
ranges are incomparable is easily seen to be equivalent to the following
property of a basis $(e_n)$, that we call {\em tight by range}. If $(y_n)$ is a
block sequence of $(e_n)$ and $A\subseteq \N$ is infinite, then
$$
[y_n]_{n\in \N}\not\sqsubseteq [e_n\del n\notin \bigcup_{i\in A}{\rm range} \;y_i].
$$
Thus, $(e_n)$ is tight by range if it is tight and for all block sequences
$(y_n)$ of $(e_n)$ the corresponding sequence of intervals $I_i$ is given by
$I_i={\rm range}\; y_i$. This property is also weaker than disjointly supported
subspaces being incomparable, which we shall call {\em tight by support}. It is
trivial to see that a basis, tight by range, is continuously tight.

We say that a basic sequence $(e_n)$ is   {\em subsequentially minimal} if any
subspace of $[e_n]$ contains an isomorphic copy of a subsequence of $(e_n)$. It
is clearly a weak form of minimality.

In \cite{KLMT} the authors study another notion in the context of certain
partly modified mixed Tsirelson spaces that they also call subsequential
minimality. According to their definition, a basis $(e_n)$ is subsequentially
minimal if any block basis has a further block basis equivalent to a
subsequence of $(e_n)$. However, in all their examples the basis $(e_n)$ is
weakly null and it is easily seen that whenever this is the case the two
definitions agree. They also define $(e_n)$ to be strongly non-subsequentially
minimal if any block basis contains a further block basis that has no further
block basis  equivalent to a subsequence of $(e_n)$. By Proposition
\ref{opposite}, this is seen to be weaker than tightness by range.

We shall now proceed to show a dichotomy between tightness by range and subsequential minimality.

\begin{thm}[4th dichotomy]\label{reflexionsdansunbus}
Let $E$ be a Banach space with a basis $(e_n)$. Then there exists a block
sequence $(x_n)$ of $(e_n)$ with one of the following properties, which are
mutually exclusive and both possible:
\begin{enumerate}
  \item Any two block subspaces of $[x_n]$ with disjoint ranges are incomparable.
  \item The basic sequence $(x_n)$ is subsequentially minimal.
\end{enumerate}
\end{thm}

Arguably Theorem \ref{reflexionsdansunbus} is not  a dichotomy in Gowers'
sense, since property (2) is not hereditary: for example the universal basis of
Pe\l czy\'nski \cite{pe} satisfies (2) while admitting subsequences with
property (1). However, it follows obviously from Theorem
\ref{reflexionsdansunbus} that any basis $(e_n)$ either has a block basis such
that any two block subspaces with disjoint ranges are incomparable or has a
block basis $(x_n)$ that is {\em hereditarily subsequentially minimal}, i.e.,
such that any block basis has a further block basis that is subsequentially
minimal. Furthermore, by an easy improvement of our proof or directly by
Gowers' second dichotomy, if the first case of Theorem
\ref{reflexionsdansunbus} fails, then one can also suppose that $[x_n]$ is
quasi minimal.

We shall call a basis $(x_n)$ {\em sequentially minimal} if  it is both
hereditarily subsequentially minimal and quasi minimal. This is  equivalent to
any block basis of $(x_n)$  having a further block basis   $(y_n)$ such
that every subspace of $[x_n]$ contains an equivalent copy of a subsequence of $(y_n)$.
We may therefore see Theorem \ref{reflexionsdansunbus} as providing a
dichotomy between tightness by range and sequential minimality.

\

Before giving the proof of Theorem \ref{reflexionsdansunbus}, we first need to
state an easy consequence of the definition of Gowers' game.

\begin{lemme}\label{independence of I}
Let $E$ be a space with a basis and assume II has a winning strategy in Gowers'
game in $E$ to play in some set $\B$. Then there is a non-empty tree $T$ of
finite block sequences such that $[T]\subseteq \B$ and for all
$(y_0,\ldots,y_m)\in T$ and all block sequences $(z_n)$ there is a block
$y_{m+1}$ of $(z_n)$ such that $(y_0,\ldots,y_m, y_{m+1})\in T$.
\end{lemme}

\begin{proof}
Suppose $\sigma$ is the strategy for II. We
define a pruned tree $T$ of  finite block bases
$(y_0,\ldots,y_m)$ and a function $\psi$ associating to each
$(y_0,\ldots,y_m)\in T$ a sequence $(z_0,\ldots,z_k)$ such that for some $k_0<\ldots<k_m=k$,
$$
\begin{array}{ccccccccccccccc}
 {\bf I} &z_0&\ldots&z_{k_0}&       &   z_{k_0+1}&\ldots&z_{k_1}&       &\ldots& & z_{k_{m-1}+1} & \ldots&
z_{k_m}&\\
 {\bf II}&      &         &              &y_0&                     &          &             &y_1&\ldots&  &&  && y_{m}\\
\end{array}
$$
has been played according to $\sigma$.

\begin{itemize}
  \item The empty sequence $\tom$ is in $T$ and $\psi(\tom)=\tom$.
  \item If $(y_0,\ldots,y_m)\in T$ and
  $$
  \psi(y_0,\ldots,y_m)=(z_0,\ldots,z_k),
  $$
  then we let
  $(y_0,\ldots,y_m,y_{m+1})\in T$ if there are some $z_k<z_{k+1}<\ldots<z_l$ and $k_0<\ldots<k_m=k$ such that
$$
\begin{array}{cccccccccccccccccc}
 {\bf I} &z_0&\ldots&z_{k_0}&       &   z_{k_0+1}&\ldots&z_{k_1}&       &\ldots& & z_{k_m+1} & \ldots&
z_{l}&\\
 {\bf II}&      &         &              &y_0&                     &          &             &y_1&\ldots&  &&  && y_{m+1}\\
\end{array}
$$
has been played according to $\sigma$   and in this case we let
  $$
  \psi(y_0,\ldots,y_m,y_{m+1})=(z_0,\ldots,z_k,z_{k+1},\ldots,z_l)
  $$
  be some such sequence.
\end{itemize}
Now, if $(y_0,y_1,y_2,\ldots)$ is such that $(y_0,\ldots,y_m)\in T$ for all
$m$, then $\psi(\tom)\subseteq \psi(y_0)\subseteq \psi(y_0,y_1)\subseteq
\ldots$ and $(y_i)$ is the play of II according to the strategy $\sigma$ in
response to $(z_i)=\bigcup_n\psi(y_0,\ldots,y_n)$ being played by I. So
$[T]\subseteq \B$. It  also follows by the construction that for each
$(y_0,\ldots,y_m)\in T$ and block sequence $(z_i)$ there is a block $y_{m+1}$
of $(z_i)$ such that $(y_0,\ldots,y_m,y_{m+1})\in T$.
\end{proof}

We now pass to the proof of Theorem \ref{reflexionsdansunbus}.

\begin{proof}
If $E$ contains $c_0$ the theorem is trivial. So suppose not. By the solution
to the distortion problem  and by passing to a subspace, we can suppose there
are two positively separated inevitable closed subsets $F_0$ and $F_1$ of the
unit sphere of $E$, i.e., such that ${\rm dist}(F_0,F_1)>0$ and every block
basis has block vectors belonging to both $F_0$ and $F_1$.

Suppose that $(e_n)$ has no block sequence satisfying (1). Then for all block
sequences $(x_n)$ there are further block sequences $(y_n)$ and $(z_n)$ with
disjoint ranges such that $[y_n]\sqsubseteq [z_n]$. We claim that there is a
block sequence $(f_n)$ and a constant $K$ such that for all block sequences
$(x_n)$ of $(f_n)$ there are further block sequences $(y_n)$ and $(y_n)$ with
disjoint ranges such that $[y_n]\sqsubseteq_K [z_n]$. If not, we can construct
a sequence of block sequences $(f_n^K)$ such that $(f_n^{K+1})$ is a block of
$(f_n^K)$ and such that any two block sequences of $(f_n^K)$ with disjoint
ranges are $K^2$-incomparable. By Lemma \ref{block uniformisation}, we then
find a block sequence $(g_n)$ of $(e_n)$ that is $\sqrt K$-equivalent with a
block sequence of $(f_n^K)$ for every $K\geq 1$. Find now block subspaces
$(y_n)$ and $(z_n)$ of $(g_n)$ with disjoint ranges and a $K$ such that
$[y_n]\sqsubseteq_{\sqrt K} [z_n]$. Then $(g_n)$ is $\sqrt K$-equivalent with a
block sequence of $(f_n^K)$ and hence we can find $K^{3/2}$-comparable block
subspaces of $(f_n^K)$ with disjoint ranges, contradicting our assumption.

So suppose $(f_n)$ and $K$ are chosen as in the claim. Then for all block
sequences $(x_n)$ of $(f_n)$ we can find an infinite set $B\subseteq \N$ and a
block sequence $(y_n)$ of $(x_n)$ such that $[y_n]_{n \in B}$ $K$-embeds into
$[y_n]_{n\notin B}$.

We claim that any block basis of $(f_n)$ has a further block basis in the
following set of normalised block bases of $(f_n)$:
$$
\A=\{(y_n)\del \a n\; y_{2n}\in F_0\cup F_1\;\&\; \e^\infty n\; y_{2n}\in
F_0\;\&\;[y_{2n+1}]_{y_{2n}\in F_0}\sqsubseteq_{K} [y_{2n+1}]_{y_{2n}\in
F_1}\}.
$$
To see this, suppose that $(x_n)$ is a block sequence of $(f_n)$ and let
$(z_n)$ be a block sequence of $(x_n)$ such that $z_{3n}\in F_0$ and
$z_{3n+1}\in F_1$. We can now find an infinite set $B\subseteq \N$ and a block
sequence $(v_n)$ of $(z_{3n+2})$ such that $[v_n]_{n \in B}\sqsubseteq_K
[v_n]_{n\notin B}$. Let now $y_{2n+1}=v_n$ and notice that we can choose
$y_{2n}=z_i\in F_0$ for $n\in B$ and $y_{2n}=z_i\in F_1$ for $n\notin B$ such
that $y_0<y_1<y_2<\ldots$. Then $(y_n)\in \A$.

Choose now a sequence $\Delta=(\delta_n)$ of positive reals, $\delta_n<{\rm
  dist}(F_0,F_1)/3$, such that if $(x_n)$ and $(y_n)$ are block bases of
$(e_n)$ with $\|x_n-y_n\|<\delta_n$, then $(x_n)\sim_2(y_n)$. Since $\A$ is
clearly analytic, it follows by Gowers' determinacy theorem that for some
block basis $(x_n)$ of $(f_n)$, II has a winning strategy to play in
$\A_\Delta$ whenever I plays a block basis of $(x_n)$. We now show that some
block basis $(v_n)$ of $(x_n)$ is
 such that any subspace of $[v_n]$ contains a sequence $2K$-equivalent to a subsequence of $(v_n)$, which will give us case (2).

Pick first by Lemma \ref{independence of I}  a non-empty tree $T$ of finite
block sequences of $(x_n)$ such that $[T]\subseteq \A_\Delta$ and for all
$(u_0,\ldots,u_m)\in T$ and all block sequences $(z_n)$ there is a block
$u_{m+1}$ of $(z_n)$ such that $(u_0,\ldots,u_m, u_{m+1})\in T$. Since $T$ is
countable, we can construct inductively a block sequence $(v_n)$ of $(x_n)$
such that for all $(u_0,\ldots,u_m)\in T$ there is some $v_n$ with
$(u_0,\ldots,u_m,v_n)\in T$.

We claim that $(v_n)$ works. For if $(z_n)$ is any block sequence of $(v_n)$,
we construct inductively a  sequence $(u_n)\in \A_\Delta$ as follows. Using
inductively the extension property of $T$, we can construct an infinite block
sequence $(h^{0}_n)$ of $(z_n)$ that belongs to $[T]$. Since  $[T]\subseteq
\A_\Delta$, there is a shortest initial segment $(u_0,\ldots, u_{2k_0})\in T$
of $(h^{0}_n)$ such that $d(u_{2k_0},F_0)<\delta_{2k_0}$. Pick now a term
$u_{2k_0+1}$ from $(v_n)$ such that $(u_0,\ldots,u_{2k_0},u_{2k_0+1})\in T$.

Again, using the extension property of $T$, there is an infinite block sequence
$(h_n^{1})$ of $(z_n)$ such that
$$
(u_0,\ldots,u_{2k_0},u_{2k_0+1})^\frown(h_n^{1})_n\in[T].
$$
 Also, as $[T]\subseteq \A_\Delta$, there is a shortest initial segment
 $$
 (u_0,\ldots,u_{2k_0},u_{2k_0+1}, \ldots,u_{2k_1})\in T
 $$
 of
 $$
 (u_0,\ldots,u_{2k_0},u_{2k_0+1})^\frown(h_n^{1})_n
 $$
that properly extends $(u_0,\ldots,u_{2k_0},u_{2k_0+1})$ and such that
$d(u_{2k_1},F_0)<\delta_{2k_1}$. We then pick a term $u_{2k_1+1}$ of $(v_n)$
such that $(u_0,\ldots,u_{2k_1},u_{2k_1+1})\in T$. We continue in the same
fashion.

At infinity, we then have a block sequence $(u_n)\in \A_\Delta$ and integers
$k_0<k_1<\ldots$ such that $d(u_{2n},F_0)<\delta_{2n}$ if and only if $n=k_i$
for some $i$ and such that for every $i$, $u_{2k_i+1}$ is a term of $(v_n)$.
Let now $(w_n)\in \A$ be such that $\|w_n-u_n\|<\delta_n$. Then, as
$\delta_n<{\rm dist}(F_0,F_1)/3$, we have that $w_{2n}\in F_0$ if and only if
$n=k_i$ for some $i$ and $w_{2n}\in F_1$ otherwise. Moreover, as $(w_n)\in \A$,
$$
[w_{2k_i+1}]_{i\in \N}=[w_{2n+1}]_{w_{2n}\in F_0}\sqsubseteq_K [w_{2n+1}]_{w_{2n}\in F_1}=[w_{2n+1}]_{n\neq k_i}.
$$
So by the choice of $\delta_n$ we have
\[\begin{split}
[u_{2k_i+1}]_{i\in \N}
&\sqsubseteq_2[w_{2k_i+1}]_{i\in \N}\sqsubseteq_K [w_{2n+1}]_{n\neq k_i}\sqsubseteq_2 [u_{2n+1}]_{n\neq k_i}.
\end{split}\]
Since $ [u_{2n+1}]_{n\neq k_i}$ is a subspace of $[z_n]$ and $(u_{2k_i+1})$ a subsequence of $(v_n)$ this finishes the proof.
\end{proof}

If, for some constant $C$, all subspaces of $[e_n]$ contain a $C$-isomorphic
copy of a subsequence of $(e_n)$, we say that $(e_n)$ is {\em subsequentially
$C$-minimal}. Our proof shows that condition (2) in Theorem
\ref{reflexionsdansunbus} may be improved to ``For some constant $C$  the
basic sequence $(x_n)$ is subsequentially $C$-minimal''.

We notice that if $(x_n)$ is hereditarily subsequentially minimal, then there
is some $C$ and a block sequence $(v_n)$ of $(x_n)$ such that $(v_n)$ is
hereditarily subsequentially $C$-minimal with the obvious definition.  To see
this, we first notice that by Proposition \ref{opposite}, $(x_n)$ can have no
block bases  $(y_n)$ such that further block subspaces with disjoint ranges are
incomparable. So, by the proof of Theorem \ref{reflexionsdansunbus}, for any
block base $(y_n)$ there is a constant $C$ and a further block basis $(z_n)$
which is subsequentially $C$-minimal. A simple diagonalisation using Lemma
\ref{block uniformisation} now shows that  by passing to a block $(v_n)$ the
$C$ can be made uniform. Recall that Gowers also proved that a quasi minimal
space must contain a further subspace which is $C$-quasi minimal
\cite{g:dicho}.

\

We also indicate a variation on Theorem \ref{reflexionsdansunbus}, relating the
Casazza property to a slightly  stronger form of sequential minimality. This
answers the original problem of Gowers left open in \cite{g:dicho}, which was
mentioned in Section \ref{intro}. This variation is probably of less interest
than Theorem \ref{reflexionsdansunbus} because the Casazza  property  does not
seem to imply tightness and also because the stronger form of sequential
minimality may look somewhat artificial (although it is satisfied by
Tsirelson's space and is reminiscent of Schlumprecht's notion of Class~1 space
\cite{S:notes}).

We say that two block sequences $(x_n)$ and $(y_n)$  {\em alternate} if either $x_0<y_0<x_1<y_1<\cdots$ or $y_0<x_0<y_1<x_1<\cdots$.
\begin{prop}

Let $E$ be a Banach space with a basis $(e_n)$. Then there exists a block
sequence $(x_n)$ with one of the following properties, which are exclusive and
both possible:
\begin{enumerate}
  \item $[x_n]$ has the Casazza property, i.e., no alternating block sequences in $[x_n]$ are equivalent.
  \item There exists a family ${\mathcal B}$ of block sequences saturating $[x_n]$
  and such that any two block sequences in ${\mathcal B}$ have subsequences
  which alternate and are equivalent.
\end{enumerate}
In particular, in case (2), $E$ contains a block subspace $U=[u_n]$ such that
for every block sequence  of $U$, there is a further block sequence equivalent
to, and alternating with, a subsequence of $(u_n)$.
\end{prop}

\begin{proof}
If $(e_n)$ does not have a block sequence satisfying (1), then any block sequence
of $(e_n)$ has a further block sequence in
$\A=\{(y_n)\del (y_{2n}) \sim (y_{2n+1})\}$. Let $\Delta$ be small enough so that
${\A}_{\Delta}=\A$. By Gowers' theorem, let $(x_n)$ be some block sequence of
$(e_n)$ so that II has a winning strategy to play in ${\A}_{\Delta}$ whenever
plays a block sequence of $(x_n)$. Let $T$ be the associated tree given by
Lemma \ref{independence of I}. By construction, for any block sequence $(z_n)$
of $(x_n)$, we may find a further block sequence $(v_n)$ such that for any
$(y_0,\ldots,y_m) \in T$, there exists some $v_n$ with
$(y_0,\ldots,y_m,v_n) \in T$. We set $f((z_n))=(v_n)$ and
${\mathcal B}=\{f((z_n))\del (z_n) \leq (x_n)\}$.
Given $(v_n)$ and $(w_n)$ in ${\mathcal B}$, it is then clear that we may find
subsequences $(v^{\prime}_n)$ and $(w^{\prime}_n)$ so that
$(v_0^{\prime},w_0^{\prime},v_1^{\prime},w_1^{\prime},\ldots) \in T$, and
therefore $(v_n^{\prime})$ and $(w_n^{\prime})$ are equivalent.
\end{proof}

We may also observe that there is no apparent relation between tightness by
range and tightness with constants. Indeed Tsirelson's space is tight with
constants and sequentially minimal. Similarly, Example \ref{gunc} is tight by
support and therefore by range, but is locally minimal. Using the techniques of
\cite{ADKM}, one can construct a space combining the two forms of tightness.

\begin{ex}\cite{exemples}
There exists a space with a basis which is tight with constants and tight by
range.
\end{ex}

Finally, if a space $X$ is locally minimal and equipped with a basis which is
tight by support and therefore unconditional (such as Example \ref{gunc}),
then the reader will easily check the following. The canonical basis of $X
\oplus X$ is tight (for a block subspace $Y=[y_n]$ of $X \oplus X$ use the sequence
of intervals associated to the ranges of $y_n$ with respect to the canonical
$2$-dimensional decomposition of $X \oplus X$), but neither tight by range nor
with constants.
However, a more interesting question remains open: does there exist a tight
space which does not contain a basic sequence which is tight by range or with
constants?

\

There is a natural strengthening of sequential minimality that has been
considered in the literature, namely, the {\em blocking principle} (also known
as the {\em shift property} in \cite{CK}) due to Casazza, Johnson, and Tzafriri
\cite{CJT}. It is known that for a normalized unconditional basis $(e_n)$ the
following properties are equivalent (see, e.g., \cite{ergodic}).
\begin{enumerate}
  \item Any block sequence $(x_n)$ spans a complemented subspace of $[e_n]$.
  \item For any block sequence $(x_n)$, $(x_{2n})\sim (x_{2n+1})$.
  \item For any block sequence $(x_n)$ and integers $k_n\in {\rm supp}\;x_n$, $(x_n)\sim (e_{k_n})$.
\end{enumerate}
Moreover, any of the above properties necessarily hold uniformly. We say that
$(e_n)$ satisfies the blocking principle if the above properties hold for
$(e_n)$. The following proposition can be proved along the lines of the proof
of the minimality of $T^*$ in \cite{CJT} (Theorem 14).

\begin{prop}
Let $(e_n)$ be an unconditional basis satisfying the blocking principle and
spanning a locally minimal space. Then $(e_n)$ spans a minimal space.
\end{prop}
Thus, by the $5$th dichotomy (Theorem \ref{5th}), we have

\begin{cor}
Let $(e_n)$ be an unconditional basis satisfying the blocking principle. Then
there is a subsequence $(f_n)$ of $(e_n)$ such that either $[f_n]$ is minimal
or $(f_n)$ is tight with constants.
\end{cor}

\section{Chains and strong antichains}\label{chains and posets}
The results in this section are in response to a question of Gowers from his
fundamental study \cite{g:dicho} and concern what types of quasi orders can be
realised as the set of (infinite-dimensional) subspaces of a fixed Banach space
under the relation of isomorphic embeddability.
\begin{prob}[Problem 7.9. in \cite{g:dicho}]\label{g:prob}
Given a Banach space $X$, let $\PP(X)$ be the set of all equivalence classes of
subspaces of $X$, partially ordered by isomorphic embeddability. For which
posets $P$ does there exist a Banach space $X$ such that every subspace $Y$ of
$X$ contains a further subspace $Z$ with $\PP(Z) = P$?
\end{prob}
Gowers  noticed himself that by a simple diagonalisation argument any such
poset $\PP(X)$ must either have a minimal element, corresponding to a minimal
space, or be uncountable. We shall now use our notion of tightness to show how
to attack this problem in a uniform way and  improve on several previous
results.

Suppose $X$ is a separable Banach space and let $SB(X)$ denote the set of all
closed linear subspaces of $X$. We equip $SB(X)$ with the so called {\em
Effros-Borel} structure, which is the $\sigma$-algebra generated by sets on the
form
$$
\{Y\in SB(X)\del Y\cap U\neq \tom\},
$$
where $U$ is an open subset of $X$. In this way, $SB(X)$ becomes a standard
Borel space, i.e., isomorphic as a measurable space to the real line equipped
with its Borel algebra. We refer to the measurable subsets of $SB(X)$ as Borel
sets. Let also $SB_\infty(X)$ be the subset of $SB(X)$ consisting of all
infinite-dimensional subspaces of $X$. Then $SB_\infty(X)$ is a Borel subset of
$SB(X)$ and hence a standard Borel space in its own right.
\begin{defi}
Suppose that $X$ is a separable Banach space and $E$ is an analytic equivalence
relation on a Polish space $\ku Z$. We say that $X$ has an $E$-{\em antichain},
if there is a Borel function $f\colon \ku Z\til SB(X)$ such that for $x,y\in\ku
Z$
\begin{enumerate}
  \item if $xEy$, then $f(x)$ and $f(y)$ are biembeddable,
  \item if $x\not\mathrel{E}y$, then $f(x)$ and $f(y)$ are incomparable.
\end{enumerate}
We say that $X$ has a {\em strong $E$-antichain} if there is a Borel function
$f\colon \ku Z\til SB(X)$ such that for $x,y\in \ku Z$
\begin{enumerate}
  \item if $xEy$, then $f(x)$ and $f(y)$ are isomorphic,
  \item if $x\not\mathrel{E}y$, then $f(x)$ and $f(y)$ are incomparable.
\end{enumerate}
\end{defi}
For example, if $=_\R$ is the equivalence relation of identity on $\R$, then
$=_\R$-antichains and strong $=_\R$-antichains simply correspond to a perfect
antichain in the usual sense, i.e., an uncountable Borel set of pairwise
incomparable subspaces. Also, having a strong $E$-antichain implies, in
particular, that $E$ Borel reduces to the isomorphism relation between the
subspaces of $X$.

The main result of \cite{flr} reformulated in this language states that if
$E_{{\fed \Sigma}_1^1}$ denotes the complete analytic equivalence relation,
then $C[0,1]$ has a strong $E_{{\fed \Sigma}_1^1}$-antichain.

We will now prove a result that simultaneously improves on two results due
respectively to the first and the second author. In \cite{ergodic}, the authors
proved that a Banach space not containing a minimal space must contain a
perfect set of non-isomorphic subspaces. This result was improved by Rosendal
in \cite{incomparable}, in which it was shown that if a space does not contain
a minimal subspace it must contain a perfect set of pairwise incomparable
spaces. And Ferenczi proved in \cite{subsurfaces} that if $X$ is a separable
space without minimal subspaces, then $E_0$ Borel reduces to the isomorphism
relation between the subspaces of $X$. Recall that  $E_0$ is the equivalence
relation defined on $2^{\N}$ by  $x E_0 y$ if and only if $\exists m\; \forall
n \geq m\; x_n=y_n$.

\begin{thm}\label{E_0-antichain}
Let $X$ be a separable Banach space. Then $X$ either contains a minimal
subspace or has a strong $E_0$-antichain.
\end{thm}

\begin{proof}Suppose $X$ has no minimal subspace.
By Theorem \ref{3rddichotomy} and Lemma \ref{ramsey tight}, we can find a basic
sequence $(e_n)$ in $X$ and a continuous function $f\colon [\N]\til [\N]$ such
that for all $A, B\in [\N]$, if $B$ is disjoint from an infinite number of
intervals $[f(A)_{2i},f(A)_{2i+1}]$, then $[e_n]_{n\in A}$ does not embed into
$[e_n]_{n\in B}$. We claim that there is a continuous function $h\colon \ca\til
[\N]$ such that
\begin{enumerate}
  \item if $xE_0y$, then $|h(x)\setminus h(y)|=|h(y)\setminus h(x)|<\infty$,
  \item if $x\not\mathrel{E_0}y$, then $[e_n]_{n\in h(x)}$ and $[e_n]_{n\in h(y)}$ are incomparable spaces.
\end{enumerate}
This will clearly finish the proof using the fact that subspaces of the same
finite codimension in a common superspace are isomorphic.

We will construct a partition of $\N$ into intervals
$$
I_0^0<I_0^1<I_0^2<I_1^0<I_1^1<I_1^2<\ldots
$$
such that if we set $J_n^0=I_n^0\cup I_n^2$ and $J_n^1=I_n^1$, the following
conditions hold:
\begin{enumerate}
  \item for all $n$, $|J_n^0|=|J_n^1|$,
  \item if $s\in 2^n$, $a=J_0^{s_0}\cup J_1^{s_1}\cup\ldots \cup J_{n-1}^{s_{n-1}}\cup I_n^0$, and $A\in [a,\N]$, then for some $i$,
$$
[f(A)_{2i},f(A)_{2i+1}]\subseteq I_n^0,
$$
  \item if $s\in 2^n$, $a=J_0^{s_0}\cup J_1^{s_1}\cup\ldots \cup J_{n-1}^{s_{n-1}}\cup I_n^1$, and $A\in [a,\N]$, then for some $i$,
$$
[f(A)_{2i},f(A)_{2i+1}]\subseteq I_n^1.
$$
\end{enumerate}
Assuming this is done, for $x\in \ca$ we set  $h(x)=J_0^{x_0}\cup
J_1^{x_1}\cup\ldots$. Then for all $n$ there is an $i$ such that
$$
[f(h(x))_{2i},f(h(x))_{2i+1}]\subseteq J_n^{x_n}.
$$

Therefore, if $x\not\mathrel{E_0}y$, then $h(y)=J_0^{y_0}\cup
J_1^{y_1}\cup\ldots$ is disjoint from an infinite number of $J^{x_n}_n$ and
thus also from an infinite number of intervals $[f(h(x))_{2i},f(h(x))_{2i+1}]$,
whence $[e_n]_{n\in h(x)}$ does not embed into $[e_n]_{n\in h(y)}$. Similarly,
$[e_n]_{n\in h(y)}$ does not embed into $[e_n]_{n\in h(x)}$.

On the other hand, if $xE_0y$, then clearly  $|h(x)\setminus
h(y)|=|h(y)\setminus h(x)|<\infty$.

It therefore only remains to construct the intervals $I_n^i$. So suppose by
induction that $I_0^0<I_0^1<I_0^2<\ldots<I_{n-1}^0<I_{n-1}^1<I_{n-1}^2$ have
been chosen (the initial step being $n=0$) such that the conditions are
satisfied. Let $m=\max J_{n-1}^0+1=\max I_{n-1}^2+1$. For each $s\in 2^n$ and
$a=J_0^{s_0}\cup J_1^{s_1}\cup\ldots \cup J_{n-1}^{s_{n-1}}$, there are by
continuity of $f$  some $k_s>m$, some interval $m\leq M_s\leq k_s$ and an
integer $i_s$ such that for all $A\in \big[a\cup\; [m,k_s],\N\big]$, we have $$
[f(A)_{2i_s},f(A)_{2i_s+1}]=M_s.
$$
Let now $k=\max_{s\in 2^n}k_s$ and $I_n^0=[m,k]$. Then if  $s\in 2^n$ and
$a=J_0^{s_0}\cup \ldots \cup J_{n-1}^{s_{n-1}}$, we have for all $A\in [a\cup
I_n^0,\N]$ some $i$ such that
$$
[f(A)_{2i},f(A)_{2i+1}]\subseteq I_n^0.
$$

Again for each $s\in 2^n$ and $a=J_0^{s_0}\cup J_1^{s_1}\cup\ldots \cup
J_{n-1}^{s_{n-1}}$ there are by continuity of $f$  some $l_s>k+1$, some
interval $k+1\leq L_s\leq l_s$ and an integer $j_s$ such that for all $A\in
\big[a\cup\; [k+1,l_s],\N\big]$, we have
$$
[f(A)_{2j_s},f(A)_{2j_s+1}]=L_s.
$$
Let now $l=\max_{s\in 2^n}l_s+k$ and $I_n^1=[k+1,l]$. Then if  $s\in 2^n$ and
$a=J_0^{s_0}\cup \ldots \cup J_{n-1}^{s_{n-1}}$, we have for all $A\in [a\cup
I_n^1,\N]$ some $j$ such that
$$
[f(A)_{2j},f(A)_{2j+1}]\subseteq I_n^1.
$$
Finally, we simply let $I_n^2=[l+1,l+|I_n^1|-|I_n^0|]$. This finishes the
construction.
\end{proof}

\begin{defi}
We define a quasi order $\subseteq^*$ and a partial order $\subseteq_0$ on the
space $[\N]$ of infinite subsets of $\N$ by the following conditions:
$$
A\subseteq^* B\equi A\setminus B \textrm{ is finite}
$$
and
$$
A\subseteq_0 B\equi \Big(A=B \textrm { or }\e n\in B\setminus A\colon\; A\subseteq B\cup [0,n[\Big).
$$
Also, if $(a_n)$ and $(b_n)$ are infinite sequences of integers, we let
$$
(a_n)\leq^*(b_n)\equi \a^\infty n\; a_n\leq b_n.
$$
\end{defi}
We notice that $\subseteq_0$ is a partial order refining the quasi order
$\subseteq^*$, namely, whenever $A\subseteq^*B$ we let $A\subseteq_0B$ if $B
\not\subseteq^* A$ or $A=B$ or  $A\triangle B$ admits a greatest element which
belongs to $B$.

\begin{prop}\label{orders}
\begin{enumerate}
  \item Any closed partial order on a Polish space Borel embeds into $\subseteq_0$.
  \item Any partial order on a set of size at most $\aleph_1$ embeds into $\subseteq_0$.
  \item The quasi order $\subseteq^*$ embeds into $\subseteq_0$, but does not Borel embed.
  \item And finally $\subseteq_0$ Borel embeds into $\subseteq^*$.
\end{enumerate}
\end{prop}

\begin{proof}

(1) By an unpublished result of A. Louveau \cite{louveau}, any closed partial
order on a Polish space Borel embeds into $(\ku P(\N),\subseteq)$. And if we
let $(J_n)$ be a partition of $\N$ into countable many infinite subsets, we see
that $(\ku P(\N),\subseteq)$ Borel embeds into $\subseteq^*$ and $\subseteq_0$
by the mapping $A\mapsto \bigcup_{n\in  A}J_n$.

(2) $\&$ (3) It is well-known that any partial order of size at most $\aleph_1$
embeds into $\subseteq^*$ and if we let $s\colon [\N]\til [\N]$ be any function
such that $|A\triangle B|<\infty\equi s(A)=s(B)$ and $|A\triangle
s(A)|<\infty$, i.e., $s$ is a selector for $E_0$, then $s$ embeds $\subseteq^*$
into $\subseteq_0$. To see that there cannot be a Borel embedding of
$\subseteq^*$ into $\subseteq_0$, we notice that if $h\colon [\N]\til [\N]$ was
a Borel function such that $A\subseteq^*B\equi h(A)\subseteq_0h(B)$,  then, in
particular, $|A\triangle B| \textrm { is finite }\equi h(A)=h(B)$,
contradicting that $E_0$ is a non-smooth equivalence relation on $[\N]$.

(4) To see that $\subseteq_0$ Borel embeds into $\subseteq^*$, we define for an
infinite subset $A$ of $\N$ a sequence of integers $g(A)=(a_n)$ by
$$
a_n=\sum_{i\in A\cap [0,n]}2^i.
$$
Suppose now that $g(A)=(a_n)$ and $g(B)=(b_n)$. Then for each $n$,
$$
a_n=b_n\equi A\cap[0,n]=B\cap [0,n]
$$
and
$$
a_n<b_n\equi \e m\in B\setminus A,\; m\leq n,\; A\cap [0,n]\subseteq B\cup [0,m[.
$$
Thus, we have $a_n=b_n$ for infinitely many $n$ if and only if $A=B$, and if
$a_n<b_n$ for infinitely many $n$, then either $B\setminus A$ is infinite or
for some $m\in B\setminus A$ we have $A\subseteq B\cup [0,m[$. Moreover, if
$B\setminus A$ is infinite, then for infinitely many $n$, $a_n<b_n$. So
$$
B\not\subseteq^*A\saa (b_n)\not\leq^*(a_n)\saa \big( B\not\subseteq^*A \textrm{ or } A\subseteq_0B\big),
$$
and thus by contraposition
$$
(b_n)\leq^*(a_n)\saa B\subseteq^*A.
$$
Also, if $(b_n)\not\leq^*(a_n)$, then $B\not\subseteq^*A$ or $A\subseteq_0B$,
so if moreover $B\subseteq_0A$,  we would have $A\subseteq_0B$ and hence $A=B$,
contradicting $g(B)=(b_n)\not\leq^*(a_n)=g(A)$. Thus,
$$
B\subseteq_0A\saa (b_n)\leq^*(a_n).
$$
To see that also
$$
(b_n)\leq^*(a_n)\saa B\subseteq_0A,
$$
notice that if $(b_n)\leq^*(a_n)$ but $B\not\subseteq_0A$, then, as
$B\subseteq^*A$, we must have $A\subseteq_0B$ and hence $(a_n)\leq^*(b_n)$. But
then $a_n=b_n$ for almost all $n$ and thus $A=B$, contradicting
$B\not\subseteq_0A$. Therefore,
$$
B\subseteq_0A\equi (b_n)\leq^*(a_n),
$$
and  we thus have a Borel embedding of $\subseteq_0$ into the quasi order
$\leq^*$ on the space $\N^\N$. It is well-known and easy to see that this
latter Borel embeds into $\subseteq^*$ and hence so does $\subseteq_0$.
\end{proof}

\begin{prop}\label{poset}
Any Banach space without a minimal subspace contains a subspace with an F.D.D.
$(F_n)$ satisfying one of the two following properties:
\begin{itemize}
  \item[(a)] if $A,B\subseteq \N$ are infinite, then
  $$
  \sum_{n\in A}F_n\sqsubseteq \sum_{n\in B}F_n\equi A\subseteq^*B,
  $$
    \item[(b)] if $A,B\subseteq \N$ are infinite, then
  $$
  \sum_{n\in A}F_n\sqsubseteq \sum_{n\in B}F_n\equi A\subseteq_0B.
  $$
\end{itemize}
\end{prop}

\begin{proof}
Suppose $X$ is a Banach space without a minimal subspace. Then by Theorem
\ref{3rddichotomy}, we can find a continuously tight basic sequence $(e_n)$ in
$X$. Using the infinite Ramsey Theorem for analytic sets, we can also find an
infinite set $D\subseteq \N$ such that
\begin{itemize}
  \item[(i)] either for all infinite $B\subseteq D$, $[e_i]_{i\in B}$ embeds into its hyperplanes,
  \item[(ii)] or for all $B\subseteq D$, $[e_i]_{i\in B}$ is not isomorphic to a proper  subspace.
\end{itemize}
And, by Lemma \ref{ramsey tight},  we can after renumbering the sequence
$(e_n)_{n\in D}$ as $(e_n)_{n\in \N}$ suppose that  there is a continuous
function $f\colon [\N]\til [\N]$ that for $A, B\in [\N]$, if $B$ is disjoint
from an infinite number of intervals $[f(A)_{2i},f(A)_{2i+1}]$, then
$[e_n]_{n\in A}$ does not embed into $[e_n]_{n\in B}$.

We now construct a partition of $\N$ into intervals
$$
I_0<I_1<I_2<\ldots
$$
such that the following conditions hold:
\begin{itemize}
  \item[-] for all $n$, $|I_0\cup\ldots \cup I_{n-1}|<|I_n|$,
  \item[-] if  $A\in [\N]$ and $I_n\subseteq A$, then for some $i$,
$$
[f(A)_{2i},f(A)_{2i+1}]\subseteq I_n.
$$
\end{itemize}
Suppose by induction that $I_0<I_1<\ldots<I_{n-1}$ have been chosen such that
the conditions are satisfied. Let $m=\max I_{n-1}+1$. For each $a\subseteq
[0,m[$ there are by continuity of $f$  some $l_a>m$, some interval $m\leq
M_a\leq l_a$ and an integer $i_a$ such that for all $A\in \big[a\cup\;
[m,l_a],\N\big]$, we have
$$
[f(A)_{2i_a},f(A)_{2i_a+1}]=M_a.
$$
Let now $l> \max_{a\subseteq [0,m[}l_a$ be such that $|I_0\cup\ldots \cup
I_{n-1}|<l-m$, and set $I_n=[m,l[$. Then if  $a\subseteq[0,m[$, we have for all
$A\in [a\cup I_n,\N]$ some $i$ such that
$$
[f(A)_{2i},f(A)_{2i+1}]\subseteq I_n,
$$
which ends the construction.

Let now $F_n=[e_i]_{i\in I_n}$. Clearly,  $\sum_{i=0}^{n-1}{\rm dim}F_i<{\rm
dim}F_n$, and if $A\setminus B$ is infinite and we let $A^*=\bigcup_{n\in A}
I_n$ and $B^*=\bigcup_{n\in B}I_n$, then $B^*$ will be disjoint from an
infinite number of the intervals defined by $f(A^*)$ and hence $\sum_{n\in
A}F_n=[e_n]_{n\in A^*}$ does not embed into $\sum_{n\in B}F_n=[e_n]_{n\in
B^*}$.

In case of (i) we have that for all infinite $C\subseteq \N$,
$$
(e_n)_{n\in C}\sqsubseteq(e_n)_{n\in C'}\sqsubseteq(e_n)_{n\in C''}\sqsubseteq(e_n)_{n\in C'''}\sqsubseteq\ldots,
$$
where $D'$ denotes $D\setminus \min D$. So, in particular, for any infinite
$A\subseteq \N$, $ \sum_{n\in A}F_n$ embeds into all of its finite
codimensional subspaces and thus if $A\setminus B$ is finite, then $\sum_{n\in
A}F_n \sqsubseteq \sum_{n\in B}F_n$. This gives us (a).

In case (ii), if $A\subseteq_0B$ but $B\not\subseteq_0 A$, we have, as ${\rm
dim}F_n>\sum_{i=0}^{n-1}{\rm dim}F_i$, that $\sum_{n\in A}F_n$ embeds as a
proper subspace of $\sum_{n\in B}F_n$. Conversely, if $\sum_{n\in
A}F_n\sqsubseteq\sum_{n\in B}F_n$, then $A\setminus B$ is finite and so either
$A\subseteq_0B$ or $B\subseteq_0A$. But if $B\subseteq_0A$ and
$A\not\subseteq_0B$, then $\sum_{n\in B}F_n$ embeds as a proper subspace into
$\sum_{n\in A}F_n$ and thus also into itself, contradicting (ii). Thus,
$A\subseteq_0B$. So assuming (ii) we have the equivalence in (b). \end{proof}

We may observe that  Tsirelson's space satisfies case (a) of Proposition
\ref{poset}, while case (b) is verified by Gowers--Maurey's space, or more
generally by any space of type (1) to (4).

By Proposition \ref{orders} and Proposition \ref{poset} we now have the
following result.

\begin{thm}
Let $X$ be an infinite-dimensional separable Banach space without a minimal
subspace and let $SB_\infty(X)$ be the standard Borel space of
infinite-dimensional subspaces of $X$ ordered by the relation $\sqsubseteq$ of
isomorphic embeddability. Then $\subseteq_0$ Borel embeds into $SB_\infty(X)$
and by consequence
\begin{itemize}
  \item[(a)] any partial order of size at most $\aleph_1$ embeds into $SB_\infty(X)$,
  \item[(b)] any closed partial order on a Polish space Borel embeds into $SB_\infty(X)$.
\end{itemize}
\end{thm}
We notice that this proves a strong dichotomy for the partial orders of Problem
\ref{g:prob}, namely, either they must be of size $1$ or must contain any
partial order of size at most $\aleph_1$ and any closed partial order on a
Polish space. In particular, in the second case we have well-ordered chains of
length $\omega_1$ and also $\R$-chains. This completes the picture of
\cite{survey}.

\section{Refining Gowers' dichotomies}\label{gowers'dichotomies}

We recall the list of inevitable classes of subspaces contained in a Banach
space given by Gowers in \cite{g:dicho}. Remember that a space is said to be
quasi minimal if any two subspaces have a common $\sqsubseteq$-minorant, and
strictly quasi minimal if it is quasi minimal but does not contain a minimal
subspace.  Also two spaces are incomparable in case neither of them embeds into
the other, and  totally incomparable if no space embeds into both of them.

\begin{thm}[Gowers \cite{g:dicho}] \label{gowers}
 Let $X$ be an infinite dimensional Banach space. Then $X$ contains a subspace $Y$ with one of the following
properties, which are all possible and mutually exclusive.
\begin{enumerate}
\item[(i)]  $Y$ is hereditarily indecomposable,
\item[(ii)] $Y$ has an unconditional basis such that any two disjointly supported block subspaces are incomparable,
\item[(iii)] $Y$ has an unconditional basis and is strictly quasi minimal,
\item[(iv)] $Y$ has an unconditional basis and is minimal.
\end{enumerate}
\end{thm}
Here the condition of (ii) that any two disjointly supported block subspaces
are incomparable, i.e., tightness by support, is equivalent to the condition
that any two such subspaces are totally incomparable or just non-isomorphic.

Theorem \ref{main} improves the list of Gowers in case (iii). Indeed, any
strictly quasi minimal space contains a tight subspace, but the space
$S(T^{(p)})$, $1<p<+\infty$ is strictly quasi minimal and not tight: it is
saturated with subspaces of $T^{(p)}$, which is strictly quasi minimal, and, as
was already observed, it is not tight because its canonical basis is symmetric.

Concerning case (i),  properties of HI spaces imply that any such space
contains a tight subspace, but it remains open whether every HI space with a
basis is tight.

\begin{quest} Is every HI space with a basis tight? \end{quest}

Using Theorem \ref{main} and Theorem \ref{main2}, we refine the list of
inevitable spaces of Gowers to 6 main classes as follows.

\begin{thm} \label{gowersbis}
Let $X$ be an infinite dimensional Banach space. Then $X$ contains a subspace
$Y$ with one of the following properties, which are all mutually exclusive.
\begin{enumerate}
\item  $Y$ is hereditarily indecomposable and has a basis such that any two block subspaces with disjoint ranges are incomparable,
\item $Y$ is hereditarily indecomposable and has a basis which is tight and sequentially minimal,
\item $Y$ has an unconditional basis such that any two disjointly supported block subspaces are incomparable,
\item $Y$ has an unconditional basis such that any two block subspaces with disjoint ranges are incomparable, and is quasi minimal,
\item $Y$ has an unconditional basis which is tight and sequentially minimal,
\item $Y$ has an unconditional basis and is minimal.
\end{enumerate}\end{thm}

We conjecture that the space of Gowers and Maurey is of type (1), although we
have no  proof of this fact. Instead, in \cite{exemples} we prove that an
asymptotically unconditional  HI space constructed by Gowers
\cite{g:asymptotic} is of type (1).

We do not know whether type (2) spaces exist. If they do, they may be thought
of as HI versions of type (5) spaces, i.e., of Tsirelson like spaces, so one
might look for an example in the family initiated by the HI asymptotically
$\ell_1$ space of Argyros and Deliyanni, whose ``ground'' space is a mixed
Tsirelson's space based on the sequence of Schreier families \cite{AD}.

The first example of type (3) was built by Gowers \cite{g:hyperplanes} and
further analysed in \cite{GM2}. Other examples are constructed in
\cite{exemples}.

Type (4) means that for any two block subspaces $Y$ and $Z$ with disjoint
ranges, $Y$ does not embed into $Z$, but some further block subspace $Y'$ of
$Y$ does ($Y'$ therefore has disjoint support but not disjoint range from $Z$).
It is unknown whether there exist spaces of type (4). Gowers sketched  the
proof of a weaker result, namely the existence of a strictly quasi minimal
space with an unconditional basis and with the Casazza property, i.e., such
that for no block sequence the sequence of odd vectors is equivalent to the
sequence of even vectors, but his example was never actually checked.
Alternatively, results of \cite{KLMT} Section 4 suggest that a mixed Tsirelson
space example might be looked for.

The main example of  a space of type (5)  is Tsirelson's space. Actually since
spaces of type (1) to (4) are either HI or satisfy the Casazza property, they
are never isomorphic to a proper subspace. Therefore, for example, spaces with
a basis saturated with block subspaces isomorphic to their hyperplanes must
contain a subspace of type (5) or (6). So our results may reinforce the idea
that Tsirelson's space is the canonical example of classical space without a
minimal subspace.

It is worth noting that as a consequence of the Theorem of James, spaces of
type (3), (4) and (5) are always reflexive.

\

Using some of the additional dichotomies, one can of course refine this picture
even further. We shall briefly consider how this can be done using the 5th
dichotomy plus a stabilisation theorem of A. Tcaciuc \cite{T} generalising a
result of \cite{FFKR}.

We state a slightly stronger version of the theorem of Tcaciuc than what is
proved in his paper and also point out that there is an unjustified use of a
recent result of Junge, Kutzarova  and Odell in his paper; their result only
holds for $1\leq p<\infty$. Tcaciuc's theorem states that any Banach space
contains either a strongly asymptotically $\ell_p$ subspace, $1 \leq p \leq
+\infty$, or a subspace $Y$ such that
$$
\a M\; \e n\; \a U_1,\ldots, U_{2n}\subseteq Y\; \e x_i\in \ku S_{U_i}\; (x_{2i-1})_{i=1}^n\not\sim_M(x_{2i})_{i=1}^n,
$$
where the $U_i$ range over infinite-dimensional subspaces of $Y$. The second
property in this dichotomy will be called {\em uniform inhomogeneity}. As
strongly asymptotically $\ell_p$ bases are unconditional, while the HI property
is equivalent to uniform inhomogeneity with $n=2$ for all $M$, Tcaciuc's
dichotomy is only relevant for spaces with an unconditional basis.

When combining Theorem \ref{gowersbis}, Tcaciuc's result, Proposition
\ref{dfko} (see also \cite{DFKO}), the 5th dichotomy, the fact that
asymptotically $\ell_{\infty}$ spaces are locally minimal, and the classical
Theorem of James, we obtain 19 inevitable classes of spaces and examples for 8
of them. The class (2) is divided into two subclasses and the class (4) into
four subclasses, which are not made explicit here for lack of an example of
space of type (2) or (4) to begin with. Recall that the spaces contained in any
of the 12 subclasses of type (1)-(4) are never isomorphic to their proper
subspaces, and in this sense these subclasses may be labeled  ``exotic''. On
the contrary ``classical'', ``pure'' spaces must belong to one of the 7
subclasses of type (5)-(6).

\begin{thm}\label{final} Any infinite dimensional Banach space contains a subspace with a basis of one of the following types:
\begin{center}
  \begin{tabular}{|l|l|l|}

    \hline

    Type       & Properties                              & Examples                                  \\

    \hline

    (1a)                     & HI, tight by range and with constants    &   ?\\

    (1b)                   & HI, tight by range, locally minimal      & $G^*$\\

    \hline

    (2) & HI, tight, sequentially minimal     &   ? \\

    \hline

    (3a) & tight by support and with constants, uniformly inhomogeneous & ? \\

 (3b) & tight by support, locally minimal, uniformly inhomogeneous & $G_u^*$ \\

 (3c) & tight by support, strongly asymptotically
                             $\ell_p$, $1 \leq p <\infty$ & $X_u$ \\

(3d) & tight by support, strongly asymptotically
                             $\ell_{\infty}$ & $X_u^*$  \\

    \hline

(4) & unconditional basis, quasi minimal, tight by range & ? \\

\hline

(5a) & unconditional basis, tight with constants, sequentially minimal,  & ? \\

     & uniformly inhomogeneous &    \\

(5b) & unconditional basis, tight, sequentially and locally minimal,
& ? \\

     & uniformly inhomogeneous &    \\

(5c) & tight with constants, sequentially minimal, & $T$, $T^{(p)}$ \\

     &  strongly asymptotically $\ell_p$, $1 \leq p<\infty$ &   \\

(5d) & tight, sequentially minimal, strongly asymptotically $\ell_{\infty}$ & ?\\

\hline

(6a) & unconditional basis, minimal, uniformly inhomogeneous & $S$ \\

(6b) & minimal, reflexive, strongly asymptotically $\ell_{\infty}$ & $T^*$\\

(6c) & isomorphic to $c_0$ or $l_p$, $1 \leq p<\infty$ & $c_0$, $\ell_p$\\

\hline
  \end{tabular}
\end{center}
\end{thm}

We know of no  space close to being  of type (5a) or (5b). A candidate for (5d)
could be the dual of some partly modified mixed Tsirelson's space not
satisfying the blocking principle (see \cite{KLMT}). Schlumprecht's space $S$
\cite{S1} does not contain an asymptotically $\ell_p$ subspace, therefore it
contains a uniformly inhomogeneous subspace, which implies by minimality that
$S$ itself is of type (6a). The definition and analysis of the spaces $G^*$,
$G^*_u$, $X_u$ and $X^*_u$ can be found in \cite{exemples}.

For completeness we should  mention that R. Wagner has also proved a dichotomy
between asymptotic unconditionality and a strong form of the HI property
\cite{W}. His result could be used to further refine the cases of type (1) and
(2).

\section{Open problems}

\begin{prob}
\begin{enumerate}
 \item Does there exist a tight Banach space admitting a basis
  which is not tight?
 \item Does there exist a tight, locally block minimal and unconditional
 basis?
\item Find a locally minimal and tight Banach space with finite
  cotype.
\item Does there exist a tight Banach space which does not  contain a basic
  sequence that is either tight by range or tight with constants? In other
  words, does there exist a locally and sequentially minimal space without a
  minimal subspace?
\item Suppose $[e_n]$ is sequentially minimal. Does there exist a block basis all of whose subsequences are subsequentially minimal?
\item Is every HI space with a basis tight?
\item Is every tight basis continuously tight?
\item Do there exist spaces of type (2), (4), (5a), (5b), (5d)?
\item Suppose $(e_n)$ is tight with constants. Does $(e_n)$ have a block sequence
that is (strongly) asymptotically $\ell_p$ for some $1\leq p< \infty$?

\item Does there exist a separable HI space $X$ such that $\subseteq^*$ Borel embeds into $SB_\infty(X)$?
\item If $X$ is a separable Banach space without a minimal subspace, does
$\subseteq^*$ Borel embed into $SB_\infty(X)$? What about more complicated
quasi orders, in particular, the complete analytic quasi order $\leq_{{\bf
\Sigma}_1^1}$?

\end{enumerate}
\end{prob}

\

\begin{flushleft}

{\em Address of V. Ferenczi:}\\
Departamento de Matem\'atica,\\
Instituto de Matem\'atica e Estat\' \i stica,\\
Universidade de S\~ao Paulo.\\
05311-970 S\~ao Paulo, SP,\\
Brazil.\\
\texttt{ferenczi@ime.usp.br}
\end{flushleft}

\

\begin{flushleft}
{\em Address of C. Rosendal:}\\
Department of Mathematics, Statistics, and Computer Science\\
University of Illinois at Chicago\\
322 Science and Engineering Offices (M/C 249)\\
851 S. Morgan Street\\
Chicago, IL 60607-7045, USA\\
\texttt{rosendal@math.uic.edu}
\end{flushleft}


\begin{thebibliography}{999}


\bibitem{T2} G. Androulakis, N. Kalton, and A. Tcaciuc, {\em On Banach spaces containing $\ell_p$ or $c_0$},
preprint.


\bibitem{AD} S. Argyros and I. Deliyanni, {\em Examples of asymptotic $\ell_1$
Banach spaces}, Trans. Amer. Math. Soc. 349 (1997), 973--995.

\bibitem{ADKM} S. Argyros, I. Deliyanni, D. Kutzarova and A. Manoussakis, {\em
Modified mixed Tsirelson spaces}, J. Funct. Anal. 159 (1998), 43--109.

\bibitem{BL} J. Bagaria and J. L\'opez-Abad,
{\em Weakly Ramsey sets in Banach spaces},  Adv. Math.  160  (2001),  no. 2, 133--174.


\bibitem{brunel}A. Brunel and L. Sucheston, {\em
On $B$-convex Banach spaces}. Math. Systems Theory 7 (1974), no. 4, 294--299.

\bibitem{CJT}P.G. Casazza, W.B. Johnson, and L. Tzafriri, {\em  On Tsirelson's space},
Israel J. Math.  47  (1984),  no. 2-3, 81--98.

\bibitem{CK} P.G. Casazza and N.J. Kalton, {\em Unconditional bases and
    unconditional finite-dimensional decompositions in Banach spaces}, Israel
    J. Math. 95 (1996), 349--373.

\bibitem{CS} P.G. Casazza and T. Shura, {\em Tsirelson's space}, Lecture Notes
  in Mathematics, 1363. Springer-Verlag, Berlin, 1989.

\bibitem{DFKO} S. Dilworth, V. Ferenczi, D. Kutzarova, and E. Odell,
{\em On strongly asymptotically $\ell_p$ spaces and minimality}, Journal of the London Math. Soc.75, 2 (2007), 409--419.



\bibitem{subsurfaces}V. Ferenczi, {\em Minimal subspaces and isomorphically homogeneous sequences in
a Banach space}, Israel J. Math.  156 (2006), 125--140.

\bibitem{exemples} V. Ferenczi and C. Rosendal, {\em Banach spaces without minimal subspaces --
Examples}, preprint.

\bibitem{flr}V. Ferenczi, A. Louveau, and C. Rosendal,
{\em The complexity of classifying separable Banach spaces up to isomorphism},
to appear in Journal of the London Mathematical Society.

\bibitem{ergodic} V. Ferenczi and C. Rosendal, {\em Ergodic Banach spaces},
Adv. Math. 195 (2005), no. 1, 259--282.

\bibitem{survey}V. Ferenczi and C. Rosendal,  {\em Complexity and homogeneity in Banach
    spaces}, Banach Spaces and their Applications in Mathematics, Ed. Beata Randrianantoanina
    and Narcisse Randrianantoanina, 2007, Walter de Gruyter, Berlin, p. 83--110.

\bibitem{FFKR}T. Figiel, R. Frankiewicz, R.A. Komorowski, and C. Ryll-Nardzewski, {\em Selecting basic
sequences in $\phi$-stable Banach spaces},
Dedicated to Professor Aleksander Pe\l czy\'nski on the occasion of his 70th birthday,
Studia Math. 159 (2003), no. 3, 499--515.

\bibitem{FLM} T. Figiel, J. Lindenstrauss and V.D. Milman, {\em The dimension
    of almost spherical sections of convex bodies}, Acta Math. 139 (1977), 53--94.

\bibitem{g:hyperplanes} W.T. Gowers, {\em A solution to Banach's hyperplane
    problem}, Bull. London Math. Soc.  26  (1994),  no. 6, 523--530.

\bibitem{g:asymptotic} W.T. Gowers, {\em A hereditarily indecomposable space with an asymptotic unconditional basis},
Geometric aspects of functional analysis (Israel, 1992--1994),  112--120, Oper. Theory Adv. Appl., 77, Birkha\"user, Basel, 1995.

\bibitem{g:hi} W.T. Gowers, {\em A new dichotomy for Banach spaces}, Geom. Funct. Anal.  6  (1996),  no. 6, 1083--1093.

\bibitem{g:dicho} W.T. Gowers, {\em An infinite Ramsey theorem and some Banach
    space dichotomies}, Ann. of Math (2) 156 (2002), 3, 797--833.

\bibitem{GM} W.T. Gowers and B. Maurey, {\em The unconditional basic sequence
    problem},  J. Amer. Math. Soc.  6  (1993),  no. 4, 851--874.



\bibitem{GM2} W. T. Gowers and B. Maurey, {\em Banach spaces with small spaces
    of operators}, Math. Ann.  307 (1997), 543--568.

\bibitem{J} W.B. Johnson, {\em A reflexive Banach space which is not
    sufficiently Euclidean}, Studia Math. 55 (1978), 201--205.

\bibitem{kechris}A. S. Kechris, {\em Classical descriptive set theory}, Springer-Verlag, New York 1995.

\bibitem{krivine} J.L. Krivine, {\em Sous-espaces de dimension finie des espaces de Banach
r\'eticul\'es}. Ann. of Math. (2)  104  (1976), no. 1, 1--29.


\bibitem{KLMT} D. Kutzarova, D. Leung, A. Manoussakis, and W.K. Tang, {\em
    Minimality properties of Tsirelson type spaces}, to appear in Studia Mathematica.

\bibitem{LT} J. Lindenstrauss and L. Tzafriri, {\em Classical Banach
spaces}, Springer-Verlag, New York, Heidelberg, Berlin (1979).


\bibitem{lopez}J. L\'opez-Abad, {\em Coding into Ramsey sets}, Math. Ann. 332 (2005), no. 4, 775--794.


\bibitem{louveau}A. Louveau, {\em Closed orders and their vicinity},  preprint 2001.

\bibitem{MMT} B. Maurey, V. Milman, and N. Tomczak-Jaegermann,
{\em Asymptotic infinite-dimensional theory of Banach spaces},
Geometric aspects of functional analysis (Israel, 1992--1994),  149--175,
Oper. Theory Adv. Appl., 77, Birkhauser, Basel, 1995.


\bibitem{OS:distortion} E. Odell and T. Schlumprecht, {\em The distortion
    problem},  Acta Math.  173  (1994),  no. 2, 259--281.

\bibitem{OS:universalunc} E. Odell and T. Schlumprecht,
{\em On the richness of the set of $p$'s in Krivine's theorem}  Geometric
aspects of functional analysis (Israel, 1992--1994),  177--198, Oper. Theory
Adv. Appl., 77, Birkh\"auser, Basel, 1995.

\bibitem{OS:universalhi} E. Odell and T. Schlumprecht, {\em A Banach space block finitely universal for monotone bases.}
Trans. Amer. Math. Soc. 352 (2000), no. 4, 1859--1888.



\bibitem{os:trees}E. Odell and T. Schlumprecht,
{\em Trees and branches in Banach spaces. } Trans. Amer. Math. Soc.  354
(2002),  no. 10, 4085--4108.



\bibitem{anna} A. M. Pe\l czar, {\em Subsymmetric sequences and minimal
spaces}, Proc. Amer. Math. Soc.  131 (2003) 3, 765-771.

\bibitem{pe} A. Pe\l czy\'nski, {\em Universal bases},  Studia Math.   32  (1969), 247--268.

\bibitem{incomparable} C. Rosendal, {\em Incomparable, non-isomorphic and minimal Banach spaces},
Fund. Math.   183  (2004) 3, 253--274.

\bibitem{asymptotic} C. Rosendal, {\em Infinite asymptotic games}, to appear in Annales de l'Institut Fourier.

\bibitem{S1} T. Schlumprecht, {\em An arbitrarily distortable Banach space},  Israel J. Math.  76  (1991),  no. 1-2, 81--95.

\bibitem{S} T. Schlumprecht, {\em How many operators exist on a Banach space?},
Trends in Banach spaces and operator theory (Memphis, TN, 2001),  295--333,
Contemp. Math., 321, Amer. Math. Soc., Providence, RI, 2003.

\bibitem{S:notes} T. Schlumprecht, unpublished notes.

\bibitem{T} A. Tcaciuc, {\em On the existence of asymptotic-$\ell_p$
    structures in Banach spaces}, Canad.  Math. Bull. 50 (2007), no. 4, 619--631.


\bibitem{tsi} B.S. Tsirelson,
{\em Not every Banach space contains $\ell_p$ or $c_0$}, Functional Anal. Appl. 8 (1974), 138--141.

\bibitem{W} R. Wagner, {\em Finite high-order games and an inductive approach towards Gowers's dichotomy},
Proceedings of the International Conference ``Analyse \& Logique'' (Mons,
1997).  Ann. Pure Appl. Logic  111  (2001),  no. 1-2, 39--60.
\end{thebibliography}
\end{document}